%% file: jw20m.tex
\documentclass[12pt,reqno,psamsfonts]{amsart}
\usepackage[alphabetic]{amsrefs}
\usepackage{feynmp}
\usepackage{latexsym}
\usepackage{amsfonts}
\usepackage{amssymb}
\usepackage{amscd}
\usepackage{amsbsy}
\usepackage{amsmath}
\usepackage{bbm}
\usepackage{euscript}
\usepackage{mathrsfs}
\usepackage{yhmath}
\usepackage[dvips]{graphicx}
\usepackage{epsfig}    
\usepackage[all]{xy}
\xyoption{knot}
\xyoption{graph}
\xyoption{tile}
\xyoption{arc}

\input ./defs/gendef.tex
\input ./defs/basdef.tex
\input ./defs/letdef.tex


\input ./defs/locdef2m.tex
\input ./defs/locdef3.tex

\input ./defs/locdefgr3.tex

\numberwithin{equation}{section}

\title[An infinite torus braid and the Jones-Wenzl projector]
{An infinite torus braid yields a categorified Jones-Wenzl projector}

\author[L.~Rozansky]{Lev Rozansky}
\address{
L.~Rozansky\\
Department of Mathematics\\
University of North Carolina at Chapel Hill\\
CB \# 3250, Phillips Hall\\
Chapel Hill, NC 27599
}
\email{rozansky@math.unc.edu}
\thanks{The work of L.R. was supported in part by the NSF grant DMS-0808974}

\begin{document}
\maketitle
\begin{abstract}
A sequence of \TLa\ elements corresponding to \cbr s with growing
twisting numbers converges to the \JWp.
We show that a sequence of
categorification complexes of these braids
also has a limit which may serve as
a categorification of the \JWp.

\end{abstract}
\tableofcontents

\section{Introduction}

A \JWp\ $\jwpn$ is a special idempotent element of the $n$-strand
\TLa\ $\cTLn$, whose defining property is
the annihilation of cap and cup tangles.
The coefficients in its expression in terms of \TLb\ tangles are
rational (rather than polynomial) functions of $q$. This suggests
that the categorification $\ctjwn$ of $\jwpn$ in the
universal tangle category $\dTLn$ constructed by D.~Bar-Natan\cx{BN1}
should be presented by a semi-infinite \chcpl. In fact, there are
two mutually dual categorifications: the complex $\ctjwpn$ which is bound from
above and the complex $\ctjwmn$ which is bound from below. We will
consider only $\ctjwpn$ in detail, since the story of $\ctjwmn$ is totally
similar.

%

The construction of $\ctjwmn$ by
B.~Cooper and S.~Krushkal\cx{CK} is based upon the
Frenkel-Khovanov formula for $\jwpn$ and requires the invention of morphisms
between constituent \TL\ tangles as well as non-trivial `thickening'
of the complex. An alternative `representation-theoretic'
approach to the categorification of the \JWp\ is developed by Igor Frenkel,
Catharina Stroppel, and Joshua Sussan\cx{FSS}.

Our approach is rather straightforward: the
categorified projector $\ctjwpn$ is a direct limit of
appropriately shifted
categorification complexes of \cbr s
(\ie braid analogs of torus links) with high \clckw\ twist (the
other projector $\ctjwmn$ comes from high \cclckw\ twists).
The limit  $\ctjwpn$
can be presented as a cone:
\xlee{eq:int1}
\ctjwpn\hteqv
\CnBv{\Ohp\big(2m(n-1)\big)\longrightarrow\cbrmns},
\xeee
where
$\gbrmn$ is a \cbr\ with $m$ full \clckw\ rotations of $n$ strands,
$\symcats{-}$ is the
categorification complex with a special grading shift, and
$\Ohp(k)$ denotes a \chcpl\ which ends at the homological degree
$-k$. Theorem\rw{th:cnpr} imposes even stronger restrictions on
the complex $\Ohp\big(2m(n-1)\big)$ in \ex{eq:int1}.


The advantage of our approach is that one can use \cbr s with high
twist as approximations to $\ctjwpn$ in a computation of \Kh\ of a
spin network which involves \JWp s:
if a spin network $\xnu$ is constructed by connecting $\jwpn$ to an \ttngnn\
$\xtau$ such that $\symcat{\xtau}\hteqv\Ohp(k)$, while a spin network $\xnum$ is constructed
by replacing $\jwpn$ in $\xnu$ with $\gbrmn$, then the homology of
$\symcat{\xnu}$ coincides with the shifted homology of
$\symcat{\xnum}$ in all homological degrees $i$
such that $i> -k - 2m(n-1)$. Thus one may say that
there is a stable limit
\xlee{eq:stlimsn}
\symcat{\xnu} =
\lim_{m\rightarrow+\infty}\symcats{\xnum}.
\xeee
We will define homological limits more precisely in
subsection\rw{sss.homcal}.

The practical
importance of the relation between $\symcat{\xnu}$
and $\symcat{\xnum}$ stems from the fact that $\xnum$ is an
ordinary link and its homology
can be computed with the help of
existing efficient computer programs even for high values
of $m$.
%
%

The simplest example of a spin network  is the unknot
`colored' by the $(n+1)$-dimensional representation of
$\mathrm{SU}(2)$ with the help of the projector $\jwpn$. Its
\Kh\ is approximated by the homology
of torus links $\mathrm{T}_{n,-mn}$ which appear as cyclic closures of
$\gbrmn$. The \Kh\ of torus links has been studied by Marko Stosic\cx{St}, who
observed that it stabilizes at lower degrees as $m$ grows. This is
a particular case of the `stable limit'\rx{eq:stlimsn}.

In Section\rw{s:notres}
we explain all notations and conventions
which are used in the paper. In particular, in
subsection\rw{sss:trgr} we define a non-traditional grading of
\Kh, which is convenient for our computations.
Then we formulate our results.


In Section\rw{s:elhomcal} we review basic facts about homological
`calculus' required to work with limits of sequences of
complexes in a homotopy category. In Section\rw{s:cbr} we construct
a sequence of categorification complexes of
\cbr s related by special \chmp s. This sequence yields $\ctjwpn$
as its direct limit. In Section\rw{s:prfs} we use homological
calculus of Section\rw{s:elhomcal} in order to prove that
$\ctjwpn$ is a categorification of the \JWp.

\def\bbS{ \mathbb{S} }
\def\So{ \bbS^1 }
\def\St{ \bbS^2 }
\def\Sot{ \So\times\St }

\subsection*{Acknowledgements}

This paper is a spinoff of a joint project with Mikhail
Khovanov\cx{KhRS} which is
dedicated to the study of categorification complexes of \cbr s and their
relation to the categorification of the Witten-Reshetikhin-Turaev
invariant of links in $\Sot$. I am deeply indebted to Mikhail for
numerous discussions and suggestions.

I would like to thank Slava Krushkal for sharing the results of
his ongoing research. I am also indebted to organizers of the M.S.R.I.
workshop `Homology Theories of Knots and Links' which stimulated
me to write this paper.

This work is supported by the NSF grant DMS-0808974.

\section{Notations and results}
\label{s:notres}
\subsection{Notations}
\label{ss:not}
\subsubsection{Tangles and \TLa}

All tangles in this paper are framed and we assume the blackboard
framing in pictures. We use the symbol 
$\xygraph{
!{0;/r1.5pc/:}
[u(0.5)]
!{\xcapv@(0)}
[u(0.45)r(0.23)]
*{\symfr\;\scriptstyle{k}}
[u(1.5)]
}
$to indicate an
addition of $k$ framing twists to a tangle strand:
\xlee{ae1.1b}
\xygraph{
!{0;/r1.5pc/:}
[u(0.5)]
!{\hover}
!{\hcap}
[u(0.5)l(0.25)]
}
\;\; = \;\;
\xygraph{
!{0;/r1.5pc/:}
[u(0.5)]
!{\xcapv@(0)}
[u(0.45)r(0.23)]
*{\symfr\;\scriptstyle{1}}
[u(1.5)]
}
\xeee

A tangle is called \emph{\plnr} if it can be presented by a diagram
without crossings. A \plnr\ tangle is called \emph{connected} or
\emph{\TLb} (\TLba) if
it does not contain disjoint circles. Let $\rTNG$ denote the set of all
framed tangles,
$\rTNGmn$ -- the set of \ttngmn s and $\rTNGn$ -- the set of
\ttngnn s.
We adopt similar notations for the set  $\rTL$ of \TLba-tangles.

We use the symbol $\tcmp$ to denote the composition of tangles:
$\xtauo\tcmp\xtaut$. The same symbol is used to denote the
multiplication in \TLa\ and the composition bifunctor in the
category $\dTL$.

A \TLa\ $\cTL$ over the ring of Laurent polynomials $\Zqqi$\footnote{It is clear from our normalization of the
Kauffman  bracket relation\rx{ae1.2} that we should rather use the
ring $\Zqqhi$. However, in all expressions in this paper the
half-integer power of $q$ appears only as a common factor, so the terms with integer
and half-integer powers of $q$ do not mix. Hence
we refer to $\Zqqi$, while keeping in mind that $\qh$ may
appear as a common factor is some expressions.}
is a quiver ring. The vertices $v_n$ of the quiver are indexed by
non-negative integers $n$ and each pair of vertices $v_m$, $v_n$,
such that $m-n$ is even, is connected
by an edge $e_{mn}$. To a vertex $v_n$ we associate a ring $\cTLnn$ (also denoted as
$\cTLn$)
and to an edge $e_{mn}$ we associate a
$\cTLn\otimes\cTLm^{\mathrm{op}}$-module $\cTLmn$. As a module,
$\cTLmn$ is generated freely by elements $\clam$ corresponding to \TL\ \ttngmn s $\xlam$, while ring
and module structures come from the composition of tangles modulo
the relation
\xlee{ae1.1}
\Bsymalg{\lcir}
 = -(\qpqi),
\xeee
which is needed to remove disjoint circles that may appear in the composition
of \TLb\ tangles.

The map $\rTNG\xrightarrow{\symalg{-}}\cTL$
associates an element $\ctau$ to a tangle $\xtau$ with the help of
\ex{ae1.1} and the Kauffman bracket relation
\xlee{ae1.2}
\Bsymalg{\xcrsp}
\;\;=\;\;
\qvh\;\;
\Bsymalg{\xpver}
\;\;+\;\;
\qvmh\;\;
\Bsymalg{\xphor}.
\xeee
This relation removes crossings and disjoint circles from the
diagram of $\xtau$, hence
\xlee{ae1.2a0}
\ctau = \sltln \xcalt\, \clam,\qquad
\xcalt = \sum_{i\in\ZZ}\xcalit\,q^i
\xeee
with only finitely many coefficients $\xcalit$ being non-zero.
%

If two tangles differ only by the framing of their strands, then
the corresponding algebra elements differ by the $q$
power factor coming from the following relation associated with
the first Reidemeister move:
\xlee{ae1.2a}
\Bsymalg{\xvfro\hspace*{-0.2cm}}
\;\; = \;\;
-q^{\frac{3}{2}}\;\;
\Bsymalg{\;\xvert\hspace*{-0.5cm}}
\xeee

A \ttngzz\ $\xL$ is a framed link, so $\symalg{\xL}$
is
the framing dependent Jones polynomial defined by the
Kauffman bracket.



We use the notations $\QcTL$ and $\cTLpinf$ for \TLa s defined over
the field $\Qq$ of rational functions of $q$ and over the field
$\Zsqqi$ of Laurent power series.
A sequence of injective homomorphisms
$\Zqqi\hookrightarrow\Qq\hookrightarrow\Zsqqi$, the latter one
generated by the expansion in powers of $q$,
produce a sequence of injective homomorphisms of the corresponding
\TLa s.

%

%

\subsubsection{The \JWp}

Let $\gcupni\in\rTLvv{n-2}{n}$ and $\gcapni\in\rTLvv{n}{n-2}$,
$1\leq i\leq n-1$, denote the following \TL\ tangles:
\ylee{ae1.3}
\gcupni=\xygraph{
!{0;/r1.5pc/:}
[r(0.25)u(0.5)]
!{\xcapv@(0)}
[u(0.5)r(1)]
*{\cdots}
[r(01)u(0.5)]
!{\xcapv@(0)}
[r(0.5)u(1)]
!{\vcap-}
[r(1.5)]
!{\xcapv@(0)}
[u(0.5)r(1)]
*{\cdots}
[r(01)u(0.5)]
!{\xcapv@(0)}
[u(1.5)l(3.5)]
*{\scriptstyle{i}}
[r(1)]
*{\scriptstyle{i+1}}
[l(3.5)]
*{\scriptstyle{1}}
[r(6)]
*{\scriptstyle{n}}
}
,
\quad\quad
\gcapni=
\xygraph{
!{0;/r1.5pc/:}
[r(0.25)u(0.5)]
!{\xcapv@(0)}
[u(0.5)r(1)]
*{\cdots}
[r(01)u(0.5)]
!{\xcapv@(0)}
[r(0.5)]
!{\vcap}
[r(1.5)u(1)]
!{\xcapv@(0)}
[u(0.5)r(1)]
*{\cdots}
[r(01)u(0.5)]
!{\xcapv@(0)}
[d(0.5)l(3.5)]
*{\scriptstyle{i}}
[r(1)]
*{\scriptstyle{i+1}}
[l(3.5)]
*{\scriptstyle{1}}
[r(6)]
*{\scriptstyle{n}}
}
\yeee
Their compositions $\xUni = \gcupni\tcmp \gcapni$ are standard
generators of the \TLa\ $\cTLn$.

The \JWp\ $\jwpn\in\QcTLn$ is the unique non-trivial idempotent element satisfying the
condition
\xlee{ae1.4}
\acapni\;\tcmp\jwpn =0,\qquad 1\leq i\leq n-1.
\xeee
The \JWp\ also satisfies the relation
\xlee{ae1.4a}
\jwpn\tcmp\;\acupni =0,\qquad 1\leq i\leq n-1.
\xeee

We denote the idempotent element of
$\cTLpinf_n$  corresponding to $\jwpn$ as $\jwpnp$.

\subsubsection{Basic notions of homological algebra}
Let $\xChA$ be a category of \chcpls\ associated with an additive
category $\xctA$. An object of $\xChA$ is a \chcpl\
\ylee{ae1.ch1}
\xbA  = (\cdots \rightarrow
\xAi\xrightarrow{\xdi}\xAio\rightarrow\cdots),
\yeee
and a morphism
between two chain complexes is a \chmp\ defined as a \mmp
\xlee{ae1.10d}
\vcenter{\xymatrix{
\xbA \ar[d]^-{\xbf} &&
\cdots\ar[r]^-{\xdimo} & \xAi \ar[r]^-{\xdi} \ar[d]^-{\yfi} &
\xAio
\ar[r]^-{\xdio} \ar[d]^{\yfio} & \cdots
\\
\xbB &&
\cdots\ar[r]^-{\xdpimo} & \xBi \ar[r]^{\xdpi} & \xBio
\ar[r]^-{\xdpio} & \cdots
}
}
\xeee
which commutes with the chain differential: $\xdpi\,\yfi = \yfio\,\xdi$ for all $i$.
The cone of a \chmp\ $\xbA\xrightarrow{\xbf}\xbB$ is a complex
\ylee{ae1.10b1}
\Cnbf
=
\lrbc{
\vcenter{
\xymatrix@C=1.5cm@R=0.5cm{
\cdots \ar[dr] \ar[r] & \xAi
\ar@{}[d] |{\oplus} \ar[r]^-{-\xdi} \ar[dr]^{-\xfi} &
\xAio
\ar@{}[d] |{\oplus}
\ar[r] \ar[dr]& \cdots
\\
\cdots \ar[r] & \xBimo\ar[r]_{\xdpimo} & \xBi \ar[r] & \cdots
}
}
}
\yeee
in which the object $\xAio\oplus\xBi$ has the homological degree
$i$.
There are two special \chmp s
$\xbB\xrightarrow{\idlbf}\Cnbf$ and
$\Cnbf[1]\xrightarrow{\chdlbf}\xbA$ associated to
the cone:
\ylee{ae1.10b2}
\xymatrix{
\xbB \ar[d]^-{\idlbf}&&
\cdots \ar[r] &
\xBi \ar[r] \ar[d]^-{0\oplus \xId}
&
\xBio \ar[r] \ar[d]^-{0\oplus \xId}
&
\cdots
\\
\Cnbf \ar[d]^-{\chdlbf} &&
\cdots \ar[r] &
\xAio \oplus \xBi \ar[r] \ar[d]^-{\xId\oplus 0} &
\xAit\oplus \xBio \ar[r] \ar[d]^-{\xId\oplus 0} &
\cdots
\\
\xbA[-1] &&
\cdots \ar[r]
&
\xAio \ar[r]
&
\xAit \ar[r]
&
\cdots
}
\yeee
These complexes and \chmp s form a \dstt:
\xlee{ae1.ch2}
\xymatrix{
\xbA\ar[r]^-{\xbf} &
\xbB \ar[r]^-{\idlbf} &
\Cnbf \ar[r]^-{\chdlbf} &
\xbA[-1]
}.
\xeee

The homotopy category of complexes $\xKhA$ has the same objects as
$\xChA$ and the morphisms are the morphisms of $\xChA$ modulo
homotopies.
%
%
We denote
homotopy equivalence by the sign $\hteqv$.
The notion of a cone extends to $\xKhA$ and there
are additional relations in that category: $\Cnv{\idlbf} \hteqv \xbA[-1]$ and
$\Cnv{\chdlbf} \hteqv \xbB[-1]$, so all vertices of a \dstt\ have
equal properties.

\subsubsection{A triply graded categorification of the Jones
polynomial}
\label{sss:trgr}
In his famous paper\cx{Kh1}, M.~Khovanov
introduced a categorification of
the Jones polynomial of links. To a diagram $\xL$ of a
link he associates a complex of graded modules
\xlee{ae1.5}
\dL = \lrbc{ \cdots \rightarrow \dLi \rightarrow \dLio\rightarrow\cdots}
\xeee
so that
if two diagrams represent the same link then the corresponding
complexes are homotopy equivalent, and the graded Euler
characteristic of $\dL$ is equal to the Jones polynomial of $\xL$.

Thus, overall, the complex\rx{ae1.5} has two gradings: the first one
is
the grading related to powers of $q$ and the second one is the
homological grading of the complex itself, the corresponding
degree being equal to $i$.
In this paper we adopt a slightly different convention which is
convenient for working with framed links and tangles. It is
inspired by matrix factorization categorification\cx{KR1} and its
advantage is that it is no longer necessary to assign orientation to
link strands in order to obtain the grading of the categorification
complex\rx{ae1.5} which would make it invariant under the second
Reidemeister move.

To a framed link
diagram $\xL$ we associate a $\ZZ \oplus\ZZ\oplus\ZZ_2$-graded complex\rx{ae1.5} with
degrees $\dgo$, $\dgt$ and $\dgh$.
The first two gradings are of the same nature as in\cx{Kh1} and, in
particular, $\dgo\dLi=i$. The third grading is an inner grading of
chain modules defined modulo 2 and of homological
nature, that is, the homological parity of an element of $\dL$,
which affects various sign factors, is the sum of $\dgo$ and
$\dgh$. Both homological degrees are either integer or
half-integer simultaneously, so the homological parity is integer
and takes values in $\ZZ_2$. The $q$-degree $\dgt$ may also take
half-integer values.


Let $\tgrshv{l}{m}{n}$ denote the shift of three degrees by $l$,
$m$ and $n$ units respectively\footnote{
Our degree shift is defined in such a way that if an object $M$
has a homogeneous degree $n$, then the shifted object $M[1]$ has a
homogeneous degree $n+1$.
}. We use abbreviated notations
$$
\tgrsshv{m}{l} = \tgrshv{m}{l}{0},\qquad
\qshv{m} = \tgrshv{m}{0}{0}
$$
as well as the following `power' notation:
$$
\tgrshv{m}{l}{n}^k = \tgrshv{km}{kl}{kn}.
$$

With new grading conventions, the categorification
formulas of\cx{Kh1} take the following form:
the module associated with an unknot is still $\ZZ[x]/(x^2)$ but with
a different degree assignment:
\begin{eqnarray}
\label{ae1.6}
&
\Bsymcat{\lcir}=\ZZ[x]/(x^2)\,
\tgrshv{-1}{0}{1},
\\
&\dgt 1 = 0, \quad \dgt x = 2,
\quad\dgo 1 = \dgo x = \dgh 1=\dgh x =0,
\end{eqnarray}
and the categorification complex of a crossing is the same as
in\cx{Kh1} but with a different degree shift:
\xlee{ae1.7}
\Bsymcat{\xcrsp}
\;\;=\;\;
\Bigg(\;\;
\Bsymcat{\xpver}
\;\tgrshv{\vthf}{-\vthf}{\vthf}
\xrightarrow{\;\;\;\;\xmrf\;\;\;\;}
\Bsymcat{\xphor}
\;\tgrshv{-\vthf}{\vthf}{-\vthf}
\vspace*{18pt}
\;\;
\Bigg),
\xeee
where $f$ is either a multiplication or a comultiplication of the
ring $\ZZ[x]/(x^2)$ depending on how the arcs in the \rhs are
closed into circles.
The resulting categorification complex\rx{ae1.5} is invariant
up to homotopy under the second and third Reidemeister moves, but
it acquires a degree shift under the first Reidemeister move:
\xlee{ae1.8}
\Bsymcat{\xvfro\hspace*{-0.2cm}}
\;\; = \;\;
\Bsymcat{\;\xvert\hspace*{-0.5cm}}
\tgrshv{\vthh}{-\vthf}{-\vthf}.
\xeee
It is easy to see that the whole categorification complex\rx{ae1.5} has a
homogeneous degree $\dgh$.

\subsubsection{A universal categorification of the \TLa}
D.~Bar-Natan\cx{BN1} described the universal category $\dTL$, whose
Grothendieck \Kzg\ is $\cTL$ considered as a $\Zqqi$-module.
We will use this category with obvious adjustments required by the new
grading conventions.

Let $\dTLt$ be an additive category whose objects are in
one-to-one correspondence with \TLb\ tangles, morphisms being generated
by tangle cobordisms (see\cx{BN1} for details). The universal category
$\dTL$ is the homotopy category of bounded complexes associated with
$\dTLt$. In other words, an object of $\dTL$ is a complex
\xlee{ae1.8a}
\xbC =
\lrbc{\cdots\rightarrow\xCi\rightarrow\xCipo\rightarrow\cdots},\qquad
\xCi =
\bigoplus_{j,\xmu}
\oltln \cjilam\,\dlam \tgrshv{j}{0}{\mu},
\xeee
where
non-negative integers $\cjilam$ are multiplicities; since the
complex is bounded, they are non-zero for only finitely many
values of $i$.

%

A categorification map $\rTNG\xrightarrow{\mpcat}\dTL$ turns a framed
tangle diagram $\xtau$ into a complex $\dtau$ according to the
rules\rx{ae1.6} and\rx{ae1.7}, the morphism $\xmrf$ in the
complex\rx{ae1.7} being the saddle cobordism. A composition of
tangles becomes a composition bi-functor
$\dTL\times\dTL\rightarrow\dTL$ if we apply
the categorified version of the rule\rx{ae1.1} in order to remove
disjoint circles:
\xlee{ae1.01}
\Bsymcat{\lcir}= \cnot \tgrshv{1}{0}{1} + \cnot\tgrshv{-1}{0}{1},
\xeee
where $\xnot$ is the empty \TL\ \ttngzz.

A complex $\dtau$ associated to a tangle $\xtau$ is defined only up to
homotopy. We use a notation $\spcc{\dtau}$ for a particular complex
with special properties which represents $\dtau$.
%

Overall,
we have the following commutative diagram:
%
\begin{equation}
\xymatrix@C=1.5cm@R=0.3cm{
& {}\dTL \ar[dd]^{\Kz}
\\
{}\rTNG \ar[ur]^{\mpcat} \ar[dr]^{\mpalg}
\\
& {}\cTL
}
\end{equation}
where the map $\Kz$ turns the complex\rx{ae1.8a} into the
sum\rx{ae1.2a0}:
\xlee{ae1.9a}
\Kz(\xbC)
=
\sltln\sum_{j}
\xcalj \,q^j\,\clam,
\qquad
\xcalj = \sum_{i,\xmu}
(-1)^{i+\xmu}\, \cjilam.
\xeee
Since the complex is bounded, the sum in the expression for
$\xcalj$ is finite.


In addition to $\dTL$ we consider
a category $\dTLp$ of complexes
bounded from above, that is, the multiplicity coefficients in the
sum\rx{ae1.8a} are zero if $i$ is greater than certain value.
Define the $q^+$ order of a \qcmd\ $\xCi$:
$\yordq{\xCi} = \xinfv{j\colon \exists\mu\colon\cjilam\neq
0}$.
A complex $\xbC$ in $\dTLp$ is \emph{\qpb} if $\lmii\yordq{\xCmi} =
+\infty$.
%
For a \qpb\ complex,
the sum in the expression\rx{ae1.9a} for $\xcalj$ is
finite, hence the element $\Kz(\xbC)$ is well defined.


\subsection{Results}
\label{ss:res}

\subsubsection{Infinite \cbr\ as a \JWp\ in a \TLa}
A braid with $n$ strands is a particular example of a \ttngnn.
A \emph{\cbr} is a braid
that can be drawn on a cylinder $\So\times[0,1]$
without intersections. In fact, all \cbr s have the form
$\btcyln^m$, $m\in\ZZ$, where $\btcyln$ is the elementary
clockwise winding \cbr:
\xlee{ae1.10p}
\btcyln \;\;=\;\;
\xygraph{
!{0;/r1.5pc/:}
!{\vover}
[u(1.5)l(0.5)]
!{\xbendr[0.5]}
[u(1.25)l(1.25)]
!{\xbendd[-0.5]}
[u(1.25)l(1)]
!{\xcapv[0.25]@(0)}
[r(1.75)]
!{\xcaph@(0)}
[u(1)]
!{\vover[-1]}
[r(1)]
!{\xbendr[-0.5]}
[u(0.75)l(0.75)]
!{\xbendd[0.5]}
[u(0.5)l(0.5)]
!{\xcapv[0.25]@(0)}
[u(1)l(1.5)]
*{\cdots}
[u(1.5)l(1.5)]
*{\cdots}
[u(0.75)l(1.5)]
*{\scriptstyle{1}}
[r(2.5)]
*{\scriptstyle{n-1}}
[r(1.25)]
*{\scriptstyle{n}}
[d(3)l(3)]
*{\scriptstyle{1}}
[r(1)]
*{\scriptstyle{2}}
[r(2.75)]
*{\scriptstyle{n}}
[u(2.2)l(0.35)]
}
\xeee
%
%
%
%
%
%
We introduce a special notation for the \cbr\ which corresponds to
$m$ full rotations of $n$ strands:
\ylee{eq:brdf}
\gbrmn = \btcyl^{mn}.
\yeee

Let $\Opqm$ denote any element of $\cTLpinf$ of the form
$\sltln\sum_{j\geq m} \xcalj\,q^j\,\clam$.
We define a \emph{\qord} of an element $\yal\in\cTLpinf$ as
$\yordq{\yal} = \xsupv{m\colon \yal = \Opqm}$.

\begin{definition}
\label{df:qlm}
A sequence of elements
\wlee{ae1.2a1}
\yal_1,\yal_2,\ldots\in\cTLpinf%
\weee
has a limit
$\lim_{k\rightarrow \infty} \yal_k = \ybet$,
if $\lmii\yordq{\ybet-\yal_k} = +\infty$.
\end{definition}

The following theorem may be known, so we do not claim
credit for it. It appears here as a by-product
and it is an easy corollary of \ex{ae2.m4}.
\begin{theorem}
\label{th:alg}
The \TL\ element corresponding to the infinite \cbr\ equals the
\JWp:
\xlee{ae1.9}
\lim_{m\rightarrow+\infty} q^{\vthf mn(n-1)}\abrmn = \jwpnp,
\xeee
where $\jwpnp\in \cTLpinf_n$ corresponds to the \JWp\
$\jwpn\in\QcTLn$.
\end{theorem}
In fact, a more general statement is also true:
\xlee{ae1.10}
\lim_{m\rightarrow+\infty} q^{\vthf m(n-1)}\symalg{\btcyln^m} =
\jwpnp,
\xeee
but its proof is more technical and we omit it here.

\subsubsection{A bit of homological calculus}
\label{sss.homcal}

Let $\xKhA$ denote the homotopy category of complexes associated
with an additive category $\xctA$ (we have in mind a particular case of $\xKhA = \dTLp$).

A \chcpl\ is considered `homologically small' if it ends at a low
(that is, high negative)
homological degree.
Let
$\Ohpm$ denote a complex which ends at $(-m)$-th homological
degree: $\Ohpm = (\cdots \xAv{-m-1} \rightarrow\xAv{-m})$. We define
a homological order of a complex $\xbA$ as $\yordh{\xbA} =
\xsupv{m\colon\xbA\hteqv\Ohpm}$.

Two complexes connected by a \chmp: $\xbA\xrightarrow{\xbf}\xbB$
are considered `homologically close' if $\Cnbf$ is homologically
small.

A \emph{\chsq} is a sequence of complexes connected by \chmp s:
\ylee{ae1.10d1}
\scA = (\xbAz\xrightarrow{\xbfz} \xbAo
\xrightarrow{\xbfo}\cdots).
\yeee
\begin{definition}
\label{df:cauchy}
A \chsq\ $\scA$ is \emph{\Cch} if $\lmii \yordh{\Cnbfi} = \pinft$.
\end{definition}
\begin{definition}
\label{df:sqlm}
A \chsq\ has a limit
\footnote{This definition differs
from the standard categorical definition of a direct limit, however
Theorem\rw{pr:spmp} indicates that our definition implies the standard one. We expect that
both definitions are equivalent.}
: $\dlm\scA = \xbA$, where $\xbA$ is a \chcpl, if
there exist \chmp s $\xbAi\xrightarrow{\xbtfi}\xbA$ such that
they form commutative triangles
\xlee{ae1.10e1}
\cmtr{\xbfi}{\xbtfio}{\xbtfi}{\xbAi}{\xbAio}{\xbA}
\xeee
%
and $\lmii\yordhr{\Cnv{\xbtfi}} =\pinft$.
\end{definition}

In Section\rw{s:elhomcal} we prove the following homology versions
of standard theorems about limits
(Propositions\rw{pr:chlm},\rw{pr:lmch} and\rw{pr:lmun}):
\begin{theorem}
\label{th:lmt}
A \chsq\ $\scA$ has a limit if and only if it is \Cch.
\end{theorem}
\begin{theorem}
\label{th:lmt2}
The limit of a \chsq\ is unique up to homotopy equivalence.
\end{theorem}

\subsubsection{Infinite \cbr\ as a \JWp\ in the universal category}

For a tangle diagram $\xtau$ let
$\dtaus$ denote the categorification complex $\dtau$
with a degree shift proportional to the number $\crsv{\xtau}$ of crossings
in the diagram $\xtau$:
\xlee{ae1.10b1}
\dtaus = \dtau\tgrshv{\vthf}{-\vthf}{\vthf}^{\crsv{\xtau}}.
\xeee
%

In subsection\rw{ss:brchsq} we define a
\chsq\ of categorification complexes of \cbr s connected by special
\chmp s
%
%
%
%
\begin{multline}
\label{ae1.10c}
\xctBn =
\Big(
\cidbrn
\xratv{\mrfz}
\cbrons \xratv{\mrfo}
\cdots
\\
\cdots
\xrightarrow{\mrfmmo}
\cbrmns \xraov{\mrfm}
\cbrmons \xrightarrow{\mrfmo}\cdots\Big).
\end{multline}
We prove that
$\yordhr{\Cnv{\mrfm}}\geq2m(n-1) + 1$,
so
$\xctBn$ is a \Csq\ and by  Theorem\rw{th:lmt} it has a unique limit:
$\dlm\xctBn =\ctjwpn \in\dTLnp $.

\begin{theorem}
\label{th:enum}
The limiting complex $\ctjwpn$ has the following properties:
\begin{enumerate}

\item A composition of $\ctjwpn$ with cap- and \uptg s is contractible:
\ylee{eq:auc}
\ccapni \;\tcmp\ctjwpn \hteqv \ctjwpn\tcmp\; \ccupni\hteqv 0.
\yeee
\item The complex $\ctjwpn$ is idempotent with respect to tangle composition:
$\ctjwpn \tcmp\ctjwpn
\hteqv \ctjwpn$.

\end{enumerate}
\end{theorem}

We provide a glimpse into the structure of $\ctjwpn$.
A complex $\xbC$  in $\cTLn$ is called \emph{\odct} if
$\gidbrn$ never appears in \qcmds\ $\xCi$.
A complex $\xbC$ in $\cTLn$ is called \emph{\otbl} if the multiplicities
$\cjilam$ of \ex{ae1.8a} satisfy the property
\xlee{ae2.m1}
\cjmilam=0\qquad\text{if $i<0$, or $j<i$, or $j>2i$.}
\xeee
%


Let $\cbrmns\xraov{\xbtfm}\ctjwpn$ be \chmp s associated
with the limit $\dlm\xctBn = \ctjwpn$ in accordance with
Definition\rw{df:sqlm}.
\begin{theorem}
\label{th:cnpr}
There exist  \odct\ \otbl\ complexes $\wbCmn$ such that
$$\Cnv{\xbtfm}\hteqv\wbCmn\spshmnm\tgrsshv{1}{-1}.$$
\end{theorem}
\noindent
In other words, there exists a distinguished triangle
%
\ylee{ae2.m2}
\wbCmn\spshmnm\qsho \xratv{\chdlbtfm} \cbrmns \xrahv{\xbtfm} \ctjwpn
\xrightarrow{\;\;\;\;\;\;\;} \wbCmn\spshmnm\tgrsshv{1}{-1}
\yeee
so there is a presentation
%
\xlee{ae2.m4}
\ctjwpn \hteqv \CnBv{ \wbCmn\spshmnm\qsho
\xrahv{\chdlbtfm} \cbrmns},
\xeee
%
where the complex $\wbCmn$ is
\odct\ and \otbl.

At $m=0$ the formula\rx{ae2.m4} becomes
\xlee{ae2.m5}
\ctjwpn \hteqv \CnBv{ \wbCzn\qsho
\xrahv{\chdlbtfz} \cidbrn},
\xeee
where the complex $\wbCzn$ is \odct\ and \otbl.

%
%

Since $\wbCzn$ is \otbl, the complex $\Cnchdlbtfz$ is
also \otbl\ and consequently \qpb. Hence $\Kctjwpn$ is well-defined. Also
$\Kctjwpn\neq 0$, because it contains $\cidbrn$ with coefficient 1.
Theorem\rw{th:enum} indicates that
$\Kctjwpn$ satisfies
defining properties of the \JWp, hence by uniqueness it is the \JWp:
\begin{corollary}
The complex $\ctjwpn$ categorifies the \JWp\ in
$\cTLpinf$:
\xlee{eq:catKz}
\Kctjwpn = \jwpn.
\xeee
\end{corollary}




\section{Elementary homological calculus}
\label{s:elhomcal}

\subsection{Limits in the category of complexes}

Consider a category $\xChA$ of \chcpls\ associated with
an additive category $\xctA$.
An $i$-th \emph{\trnc} $\xtrniv{\xbA}$ of a \chcpl\ $\xbA$ is
the \chcpl\
$\xAmi\xrightarrow{\xdmi}\xAmio\rightarrow\cdots$.
An
\trnci\ of a \chmp\ $\xbf$ is defined similarly.

Define an \emph{\isor}
$\ysiobf$
of a chain map $\xbA\xrightarrow{\xbf}\xbB$  as
the largest number $i$ for which a truncated \chmp\ $\xtrniv{\xbf}$ is an
isomorphism of truncated complexes.

\begin{remark}
\label{rm:cnord}
Consider a \dstt\rx{ae1.ch2}.
If $\xbA\hteqv\Ohpm$, then $\ysiov{\idlbf}\geq m-1$.
\end{remark}

\begin{definition}
A \chsq\
$\scA = (\xbAo\xrightarrow{\xbfo}\xbAt\xrightarrow{\xbft}\cdots)$ is
\emph{\stblz} if
$\lim_{i\rightarrow\infty} \ysiobfi=\pinft$.
\end{definition}

\begin{definition}
A \chsq\ $\scA$ has a \tchlm\ $\chlm\scA=\xbA$ if there exist
\chmp s $\xbAi\xrightarrow{\xbtfi}\xbA$ such that
$\xbtfi = \xbtfio\,\xbfi$ and $\lmii \ysiorv{\xbtfi} = \pinft$.
\end{definition}

The following two theorems are easy to prove:
\begin{theorem}
A \chsq\ has a \tchlm\ if and only if it is \stblz.
If a \tchlm\ exists then it is unique.
\end{theorem}

%

\begin{theorem}
\label{pr:fnctchlm}
Suppose that $\chlm\scA = \xbA$.
Then for a complex $\xbB$ and \chmp s
$\xbAi\xrightarrow{\ybgi}\xbB$ such that $\ybgi = \ybgio\xbfi$,
 there exists a unique \chmp\ $\xbA\xrightarrow{\ybg}\xbB$
such that $\ybgi = \ybg\,\xbtfi$.
\end{theorem}

\begin{definition}
\label{df:chlmmp}
A sequence of \chmp s $\xbA\xrightarrow{\xbfz,\xbfo,\cdots}\xbB$
has a \tchlm\ $\lmii\xbfi = \xbf$ if for any $N$ there exists $N\p$
such that $\xtrnNv{\xbfi} = \xtrnNv{\xbf}$ for any $i\geq N\p$.
\end{definition}

\subsection{Limits in the homotopy category}

Definitions\rw{df:cauchy} and\rw{df:sqlm} extend the notion of a
\stblz\ \chsq\ and its limit to the homotopy category $\xKhA$:
obviously, a \stblz\ \chsq\ is \Cch, while $\chlm\scA=\xbA$ implies
$\dlm\scA=\xbA$.

\begin{proposition}
\label{pr:chlm}
A \Csq\ has a limit.
\end{proposition}
\proof
Consider a \Csq\ $\scA$. We construct a special \chcpl\
$\xbAs$ such that $\dlm\scA=\xbAs$ in accordance with
Definition\rw{df:sqlm}. Roughly speaking, we take $\xbAz$ and
attach to it the cones $\Cnbfi$ represented by homologically small
complexes, one by one. The result is a sequence $\scAs=\xbApz,\xbApo,\ldots$ of \stblz\
complexes $\xbApi$ such that $\xbApi\hteqv\xbAi$, and
$\xbAs=\chlm\scAs$ is their \tchlm.

Here is a detailed explanation.
By Definition\rw{df:cauchy}, there exist complexes $\ybCi$ such
that
\xlee{ae1.10a1}
\Cnv{\xbfi} \hteqv \ybCi[1] = \Ohpmi,\qquad\lmii m_i=+\infty.
\xeee
The complexes $\xbAi$, $\xbAio$ and $\ybCi$ form exact triangles:
\ylee{ae1.10a2}
\xymatrix{
\ybCi \ar[r]^-{\chdlbfi} & \xbAi \ar[r]^-{\xbfi} & \xbAio \ar[r] &
\ybCi[-1]
}
\yeee
and $\xbAio \hteqv \Cnv{\chdlbfi}$. We define recursively a new sequence
of complexes $\scAp = (\xbApz \xrightarrow{\idlbgz} \xbApo\xrightarrow{\idlbgo}\cdots)$ by
the relations $\xbApz = \xbAz$, $\xbApi\hteqv\xbAi$ and
$\xbApio = \Cnv{\ybgi}$, where the \chmp\
$\ybCi\xrightarrow{\ybgi}\xbApi$ is homotopy equivalent to the
\chmp\ $\chdlbfi$. In other words,
\xlee{ae1.10g1}
\xbApio = \Cnv{\ybCi\xrightarrow{\ybgi}
\underbrace{
\Cnv{\ybCimo\xrightarrow{\ybgimo}\cdots\xrightarrow{\ybgt}
\underbrace{
\Cnv{\ybCo\xrightarrow{\ybgo}
\underbrace{
\Cnv{\ybCz\xrightarrow{\chdlbfz}\xbAz
}}_{\xbApo}\;
}
}_{\xbApt}\;
}
}_{\xbApi}\;
}
\xeee

According to Remark\rw{rm:cnord}, $\ysiov{\idlbgi}\geq m_i$,
hence the sequence
$\scAp$ is \stblz, so there exists a chain limit
$\chlm\scAp = \xbAs$ and consequently there is a limit
$\dlm\scA=\xbAs$.\qed

Simply saying, the complex $\xbAs$ is an infinite \mtcn\ extension
of the complex\rx{ae1.10g1}:
\xlee{ae1.10g2}
\xbAs =
\cdots\xrightarrow{\ybgh}\Cnv{\ybCt
\xrightarrow{\ybgt}
\Cnv{\ybCo\xrightarrow{\ybgo}
\Cnv{\ybCz\xrightarrow{\chdlbfz}\xbAz
}
}
}.
\xeee

For our applications it is important to express $\Cnbtfz$ in terms
of complexes $\ybCi$. This can be done by rearranging
the infinite \mtcn\rx{ae1.10g2} with the help of associativity of
cone formation, which exists even within the category $\xChA$:
%
%
\xlee{ae1.10h1}
\xbAs = \Cnv{\ybtC\xrightarrow{\ybtg}\xbAz},\qquad
\ybtC =\cdots\xrightarrow{\ybht}
\Cnv{\ybCt[1]\xrightarrow{\ybho}\Cnv{\ybCo[1]\xrightarrow{\ybhz}\ybCz}},
\xeee
so that $\xbtfz \hteqv \idlv{\ybtg}$, and $\Cnbtfz\hteqv\ybtC[-1]$
is expressed in terms of complexes $\ybCi$
arranged into an infinite \mtcn\ $\ybtC$. Here is a more formal
statement.
\begin{theorem}
\label{th:rshfl}
For a \Csq\ $\scA$ there exists another \Csq\
$\yctC = (\ybCz
\xrightarrow{\ybhpz} \ybtCo\xrightarrow{\ybhpo}\cdots)$
and \chmp s
$\ybCi[1]\xrightarrow{\ybhi} \ybtCi$ such that
$\Cnv{\ybhi}=\ybtCio$, $\ybhpi = \idlv{\ybhi}$ and for the
limiting complex $\ybtC=\chlm\yctC$ there exists a \chmp\
$\ybtC\xrightarrow{\ybtg}\xbAz$ such that $\xbAs = \Cnv{\ybtg}$,
$\;\xbtfz \hteqv \idlv{\ybtg}$ and consequently $\Cnbtfz\hteqv\ybtC[-1]$.
\end{theorem}

\proof
%
%
Let us recall the associativity of cones in a general setting.
For a \chmp\ $\xbA\xrightarrow{\xbf}\xbB$,
a \chmp\ $\xbC\xrightarrow{\ybg} \Cnbf$
is a sum: $\ybg = \ybgA \oplus\ybgB$
\ylee{ae1.10f1}
\xymatrix{
 & \xbA\ar[d]^-{\xbf}
\\
\xbC \ar[ur]|{[1]}^-{\ybgA} \ar[r]_-{\ybgB}
& \xbB
}
\yeee
where $\xbC\xrightarrow{\ybgA}\xbA[-1]$ is a \chmp\ and
$\xbC\xrightarrow{\ybgB}\xbB$ is a \mmp. Now it is obvious
that
\xlee{ae1.10f2}
\Cnv{\xbC\xrightarrow{\ybg}\Cnv{\xbA\xrightarrow{\xbf}\xbB}}
=\Cnv{\Cnv{\xbC[1]\xrightarrow{\ybgA}\xbA}\xrightarrow{\ybgB\oplus\xbf}\xbB
}.
\xeee

We apply the associativity relation\rx{ae1.10f2} to
\mtcn s\rx{ae1.10g1} consecutively for $i=1,2,\ldots$ in order to rearrange
them, so that $\xbApi = \Cnv{\ybtCi\xrightarrow{\ybtgi}\xbAz}$,
while the complexes $\ybtCi$ and \chmp s $\ybtgi$ are defined
recursively: $\ybtCz=\ybCz$, $\ybtgz = \chdlbfz$, $\ybtCio = \Cnv{\ybhi}$,
while the \chmp s $\ybCi[1]\xrightarrow{\ybhi} \ybtCi$ and
$\ybtCio\xrightarrow{\ybtgio}\xbAz$
are defined by applying the associativity
relation\rx{ae1.10f2} to the  double cone on the second line of
the following equation:
\begin{equation}
\label{ae1.10f3}
\begin{split}
\xbApio & = \Cnv{\ybCi\xrightarrow{\ybgi}\xbApi}
\\
& = \Cnv{\ybCi\xrightarrow{\ybgi}\Cnv{\ybtCi\xrightarrow{\ybtgi}\xbAz} }
\\
& = \Cnv{\Cnv{\ybCi[1]\xrightarrow{\ybhi} \ybtCi
} \xrightarrow{\ybtgio}
\xbAz
}
\\
& = \Cnv{\ybtCio\xrightarrow{\ybtgio}\xbAz}.
\end{split}
\end{equation}
Distinguished triangles
\ylee{ae1.10f3}
\xymatrix{
\ybCi[1]\ar[r]^-{\ybhi}
&
\ybtCi \ar[r]^-{\idlbhi}
&
\ybtCio \ar[r]
&
\ybCi
}
\yeee
determine  \chmp s $\ybhpi=\idlbhi$ of the \chsq\  $\yctC = (\ybtCz
\xrightarrow{\ybhpz} \ybtCo\xrightarrow{\ybhpo}\cdots)$. By
Remark\rw{rm:cnord} it has a \tchlm: $\chlm\yctC =
\ybtC$, which is an infinite \mtcn:
\ylee{ae1.10f4}
\ybtC =\cdots\xrightarrow{\ybht}
\Cnv{\ybCt[1]\xrightarrow{\ybho}\Cnv{\ybCo[1]\xrightarrow{\ybhz}\ybCz}}.
\yeee
The \chmp s $\ybtCi\xrightarrow{\ybhpi}\ybtCio$ satisfy a relation
$\ybtgi = \ybtgio\,\ybhpi$, so by Theorem\rw{pr:fnctchlm} there
exists a unique \chmp\ $\ybtC \xrightarrow{\ybtg}\xbAz$ such that
$\ybtgi=\ybtg\,\ybhtpi$.
It is easy to show that $\xbAs = \Cnv{\ybtC\xrightarrow{\ybtg}\xbAz}$,
and $\xbtfz = \idlv{\ybtg}$, hence $\Cnbtfz\hteqv \ybtC$.\qed

It is easy to prove the analog of Theorem\rw{pr:fnctchlm}:
\begin{theorem}
\label{pr:spmp}
For a complex $\xbB$ and \chmp s
$\xbAi\xrightarrow{\ybgi}\xbB$ such that $\ybgi \hteqv \ybgio\xbfi$,
 there exists a unique (up to homotopy) \chmp\ $\xbAs\xrightarrow{\ybg}\xbB$
which forms commutative triangles
\xlee{ae1.10f1}
\cmtr{\xbtfi}{\ybg}{\ybgi}{\xbAi}{\xbAio}{\xbB}
\xeee
%
\end{theorem}

In order to complete the proof
of Theorems\rw{th:lmt} and\rw{th:lmt2},
we need two simple propositions. The first one establishes a
triangle inequality for homological orders of cones.
\begin{proposition}
If three \chmp s form a commutative triangle
\xlee{ae1.10c1}
\xymatrix{
\xbA \ar[r]_-{\xbfAB} \ar@/^1pc/[rr]^-{\xbfAC} &
\xbB \ar[r]_-{\xbfBC} &
\xbC
},\qquad \xbfAC\hteqv \xbfBC\xbfAB.
\xeee
then the homological orders of their cones satisfy the inequalities
\begin{align}
\label{eq:in1}
\yordch{\xbfAB}\geq \min\big(\yordch{\xbfAC},\yordch{\xbfBC}-1\big),
\\
\label{eq:in2}
\yordch{\xbfBC}\geq \min\big(\yordch{\xbfAB}+1,\yordch{\xbfAC}\big).
\end{align}
\end{proposition}
\proof
If \chmp s form a commutative triangle\rx{ae1.10c1}, then their cones form a \dstt
\ylee{ae1.10c2}
\Cnv{\xbfAB}\xrightarrow{\ybgo} \Cnv{\xbfAC} \xrightarrow{\ybgt} \Cnv{\xbfBC}
\xrightarrow{\ybgh}\Cnv{\xbfAB}[1],
\yeee
so the first inequality follows from the relation
$\Cnv{\xbfAB}\hteqv\Cnv{\ybgt}[1]$ and the second inequality
follows from the relation $\Cnv{\xbfBC} \hteqv \Cnv{\ybgo}$.\qed

The second proposition says that if a complex is homologically
infinitely small then it is contractible.
\begin{proposition}
\label{pr:ismc}
If $\yordh{\xbA} = +\infty$ then $\xbA$ is contractible.
\end{proposition}
\proof Since $\yordh{\xbA} = +\infty$, there exist complexes
$\xbAi\hteqv\xbA$, such that $\xbAi=\Ohpmi$ and $\lmii m_i =
+\infty$. Consider a sequence of \chmp s establishing homotopy
equivalence between the complexes:
\ylee{ae1.10c3}
\xymatrix@C=0.5cm{
\xbA \ar@<0.6ex>[r]^-{\xbfz}
&
\xbAo
\ar@<0.3ex>[l]^-{\ybgz}
\ar@<0.6ex>[r]^-{\xbfo}
&
\xbAt
\ar@<0.3ex>[l]^-{\ybgo}
\ar@<0.6ex>[r]
&
\cdots
\ar@<0.3ex>[l]
\ar@<0.6ex>[r]
&
\xbAi
\ar@<0.3ex>[l]
\ar@<0.6ex>[r]^-{\xbfi}
&
\xbAio
\ar@<0.3ex>[l]^-{\ybgi}
\ar@<0.6ex>[r]
&
\cdots
\ar@<0.3ex>[l]
}
,\qquad
\yIdAi-\ybgi\xbfi =
\atcmv{\xbdi}{\xbhi},
\yeee
where $\yIdAi$ is the identity
\chmp\ of $\xbAi$, while $\xbAi[1]\xrightarrow{\xbhi}\xbAi$ is a homotopy \chmp\ (it does
not commute with the chain differential $\xbdi$ in the complex $\xbAi$).

Consider the compositions $\xbhfi = \xbfi\cdots\xbfo\xbfz$,
$\xbhgi = \ybgz\ybgo\cdots\ybgi$ and $\xbhhi =
\xbhgimo\,\xbhi\,\xbhfimo$. It is easy to see
that $\xbhgimo\,\xbhfimo - \xbhgi\,\xbhfi = \atcmv{\xbd}{\xbhhi}$,
hence $\yIdA - \xbhgi\,\xbhfi = \atcmv{\xbd}{\xbchi}$, where
$\xbchi = \xbhhz + \xbhho +\cdots + \xbhhi$.
There is a limit (\cf Definition\rw{df:chlmmp}) $\lmii\xbchi = \xbch$, while
$\lmii\xbhgi\,\xbhfi = 0$, hence $\yIdA = \atcmv{\xbd}{\xbch}$
which means that the complex $\xbA$ is contractible.
\qed

\begin{proposition}
\label{pr:lmch}
If a \chsq\ $\scA$ has a limit, then it is \Cch.
\end{proposition}
\proof
The inequality\rx{eq:in1} applied to the commutative
triangle\rx{ae1.10e1} says that
$$\yordch{\xbfi}\geq
\min\lrbc{\zordch{\xbtfi},\zordch{\xbtfio}-1},$$
hence the limit $\lmii\zordch{\xbtfi} = +\infty$ implies the \Cch\
property of $\scA$.

\begin{proposition}
\label{pr:lmun}
If a \chsq\ $\scA$ has a limit then it is unique.
\end{proposition}
\proof
If $\scA$ has a limit then by Proposition\rw{pr:lmch} it is \Cch.
Hence it has a special limit $\xbAs$ described in the
proof of Proposition\rw{pr:chlm}. If $\scA$ has another limit
$\xbAp$ with \chmp s $\xbAi\xrightarrow{\xbtfpi}\xbAp$ then by
Theorem\rw{pr:spmp} there is a \chmp\
$\xbAs\xrightarrow{\ybg}\xbAp$ with commutative
triangles\rx{ae1.10f1}. The inequality\rx{eq:in2} says
\ylee{ae1.10c3}
\yordch{\ybg}
\geq \min\big(\zordch{\xbtfi}+1,\zordch{\ybgi}\big).
\yeee
Since both cones in the \rhs become homologically infinitely small
at $i\rightarrow +\infty$, the cone $\Cnv{\ybg}$ is also
homologically infinitely small. Then Proposition\rw{pr:ismc} says
that $\Cnv{\ybg}$ is contractible and as a result
$\xbAp\hteqv\xbAs$.\qed

We end this section with a theorem which follows easily from
Definition\rw{df:sqlm}.
\begin{theorem}
\label{th:zl}
If a \chsq\ $\scA$ satisfies the property $\lmii\yordh{\xbAi} =
+\infty$ then its limit is contractible: $\dlm\scA = 0$.
\end{theorem}

\section{A \chsq\ of
categorification complexes of  \cbr s} 
\label{s:cbr}
\subsection{A special categorification complex of a \ngbr}

Let $\xsgi$ denote an elementary negative $n$-strand braid:
\ylee{ae2.1}
\xsgi=\xygraph{
!{0;/r1.5pc/:}
[r(0.25)u(0.5)]
!{\xcapv@(0)}
[u(0.5)r(1)]
*{\cdots}
[r(01)u(0.5)]
!{\xcapv@(0)}
[r(0.5)u(1)]
!{\vcross}
[r(1.5)u(1)]
!{\xcapv@(0)}
[u(0.5)r(1)]
*{\cdots}
[r(01)u(0.5)]
!{\xcapv@(0)}
[d(0.5)l(3.5)]
*{\scriptstyle{i}}
[r(1)]
*{\scriptstyle{i+1}}
[l(3.5)]
*{\scriptstyle{1}}
[r(6)]
*{\scriptstyle{n}}
}
\yeee

\begin{theorem}
\label{th:prop}
If an $n$-strand braid $\brb$ can be presented as a product of elementary
negative braids: $\brb = \xsgiv{k}\cdots\xsgiv{2}\xsgiv{1}$, then
its categorification complex has a special presentation $\cbrbs$:
\xlee{aea2.1}
\cbrbas =
\Big(\ldots\rightarrow\xCmt\rightarrow\xCmo\rightarrow\cidbrn\Big)
\xeee
such that the complex
\xlee{aea2.2}
\xbC = (\ldots\rightarrow\xCmt\rightarrow\xCmo)\tgrsshv{-1}{1}
\xeee
is
\odct\ and \otbl.



\end{theorem}

More abstractly, the theorem says that there exists a
\odct\ and \otbl\ complex $\xbC$ and a \chmp\
$\xbC\rightarrow\cidbrn$ such that $\cbrba \hteqv
\CnBv{\xbC\qsho\rightarrow\cidbrn}$.

\begin{remark}
\label{rm:spbl}
Theorem\rw{th:prop} implies that the special complex
$\cbrbas$ is  \otbl.
\end{remark}

\pr{Theorem}{th:prop}
Let $\xlam$ be a \TL\ \ttngnn. Fix $i$ such that $1\leq i\leq n-1$. If the
composition $\gcapni\tcmp\xlam$ does not contain a disjoint circle,
then, in accordance with \ex{ae1.7},
we define the special categorification complex of $\xsgi\tcmp\xlam$ as
%
\xlee{ae2.3}
\symcatps{\xsgi\tcmp\xlam}  =
\Big(\symbcat{\xUni\tcmp\xlam}\tgrshv{1}{-1}{1}
\rightarrow \dlam \Big)
\xeee
%
If $\gcapni\tcmp\xlam$ contains a disjoint circle, then $\xlam$ must
have the form $\gcupni\tcmp\xlamp$. Hence
$\xsgi\tcmp\xlam=\xsgi\tcmp\gcupni\tcmp\xlamp$. The tangle $\xsgi\tcmp\gcupni$
is the same as $\gcupni$ with a positive framing twist, so
according to \ex{ae1.8},
$\bsymcat{\xsgi\tcmp\gcupni} = \ccupni \tgrshv{\vthh}{-\vthf}{-\vthf}$.
Hence in this case we define the special categorification complex
of $\xsgi\tcmp\xlam$ simply as shifted $\dlam$:
\xlee{ae2.4}
\symcatps{\xsgi\tcmp\xlam} = 
\dlam
\tgrshv{2}{-1}{0}.
\xeee
%

Now we define a recursive algorithm for constructing the complex
$\cbrbas$. For $\brb = \gidbrn$ we define $\cbrbas = \cidbrn$. Let
$\brb = \xsgiv{k}\tcmp\cdots\tcmp\xsgiv{1}$ and
suppose that we have defined its special complex $\cbrbas$. We
define the special categorification complex of a
braid $\brbp=\xsgikpo\tcmp\brb$ by applying the rules\rx{ae2.3}
and\rx{ae2.4} to all constituent tangles $\xlam$ in the complex
$\cbrbs$ (see the formula\rx{ae1.8a}).

We prove the properties of $\cbrbas$ by induction over $k$.
If $k=0$ then $\brb = \gidbrn$ and the properties of $\cbrbas$ are
obvious.

Suppose that the special categorification complex
$\cbrbas$ of a braid $\brb = \xsgiv{k}\tcmp\cdots\tcmp\xsgiv{1}$
has the form\rx{aea2.1} and its tail\rx{aea2.2} is \odct\ and
\otbl.
Consider
a longer braid
$\brbp=\xsgikpo\tcmp\brb$. The object $\cidbrn$ may appear in $\cbrbpas$ if
and only if $\xlam=\gidbrn$ and the extra crossing $\xsgikpo$ is
\nsplcd\ in \ex{ae2.3}, hence
$\cbrbpas$ has the form\rx{aea2.1} and its tail\rx{aea2.2} is
\odct.

If the negative crossing $\xsgikpo$ is composed with the head
$\cidbrn$ of the complex\rx{aea2.1}, then the formula\rx{ae2.3}
applies and the tangle $\xUv{n}{i_k+1}$ appearing in the tail of
$\cbrbpas$ satisfies the property\rx{ae2.m1}.


If the crossing $\xsgikpo$ is composed with a
\TL\ tangle $\xlam$ from the $(-i)$-th \qcmd\ $\xCmi$ (see \ex{ae1.8a})
in the tail of the complex $\cbrbas$ with the $q$-degree shift $j$
satisfying the inequality $i-1 \leq j-1 \leq 2(i-1)$, then the
shifted objects in the \rhs of \eex{ae2.3} and\rx{ae2.4}
also satisfy this inequality.\qed

The picture\rx{ae1.10p} presents a \cbr\ as a product of negative
crossings, hence
\begin{corollary}
\label{cr:otbl}
A \cbr\ $\gbrmn$ has a special \otbl\ categorification complex
$\cbrmnps$. In particular, for $m=1$
\xlee{ae2.5}
\cbronps
 = \CnBv{\xbCon\qsho\rightarrow\cidbrn},
\xeee
where the complex
$\xbCon$ is \odct\ and \otbl.
\end{corollary}

\subsection{Special morphisms between \cbr\ complexes}
\label{ss:brchsq}

Relation\rx{ae2.5}  indicates that there is a \dstt\
%
%
%
\ylee{ae2.6}
\xbCon\qsho \longrightarrow
\cidbrn \xratv{\mrfo}
\cbrons \longrightarrow
\xbCon\tgrsshomo
\yeee
and
\xlee{ae2.6a}
\Cnv{\mrfo} \hteqv \xbCon\tgrsshomo.
\xeee
%
%
Composing both sides of the morphism $\mrfo$ with
the \cbr\ complex
$\cbrmns$,
we get a morphism
%
\ylee{ae2.7}
\cbrmns \xratv{\mrfm}\cbrmons
\yeee
such that
\xlee{ae2.8}
\Cnv{\mrfm} \hteqv \Cnv{\mrfo}\tcmp\cbrmns.
\xeee
\begin{theorem}
\label{th:2.1}
The cone\rx{ae2.8} can be presented by a shifted complex
\ylee{ae2.9}
\Cnv{\mrfm} \hteqv \xbCmn \tgrsshnontm\tgrsshomo,
\yeee
%
such that $\xbCmn$ is \odct\ and \otbl.
\end{theorem}

The proof is based on a simple geometric lemma:
\begin{lemma}
\label{l:1}
For $n\geq 2$, the following two compositions of framed tangles are isotopic:
%
\xlee{ae2.b}
\gcapni \tcmp\;\gbron = \gbronmt\;\tcmp\gcapnit
\xeee
where $\gcapnit$ is the tangle $\gcapni$ with double framing twist
on the cap:
\ylee{ae2.10}
\gcapnik=
\xygraph{
!{0;/r1.5pc/:}
[r(0.25)u(0.5)]
!{\xcapv@(0)}
[u(0.5)r(1)]
*{\cdots}
[r(01)u(0.5)]
!{\xcapv@(0)}
[r(0.5)]
!{\vcap}
[r(1.5)u(1)]
!{\xcapv@(0)}
[u(0.5)r(1)]
*{\cdots}
[r(01)u(0.5)]
!{\xcapv@(0)}
[d(0.5)l(3.5)]
*{\scriptstyle{i}}
[r(1)]
*{\scriptstyle{i+1}}
[l(3.5)]
*{\scriptstyle{1}}
[r(6)]
*{\scriptstyle{n}}
[l(3)u(1)]
*{\symfr}
[u(0.5)]
*{\scriptstyle{k}}
}
\yeee
\end{lemma}
\proof
This lemma is geometrically obvious: a cap on a pair of adjacent strands slides down
through the \cbr\ to the
bottom.\qed

An immediate corollary of \ex{ae2.b} and of the framing change
formula\rx{ae1.8} is the following relation:
%
\xlee{ae2.11}
\bsymcats{\gcapni \tcmp\gbrmn} \hteqv \bsymcats{\gbrmnmt\;\tcmp\gcapni}
\tgrsshnontm.
\xeee


In order to prove Theorem\rw{th:2.1}, we need three simple
propositions.
For a positive integer $d\leq \frac{n}{2}$,
let $\stI=(i_1,\ldots,i_d)$ be a sequence of positive integer
numbers such that $i_k<n-2k+2$ for all $k\in\{1,\ldots,d\}$.
A \emph{\aptg} $\gcapnI$
is a $(n,n-2d)$-tangle which
can be presented as a product of $d$ tangles of the form
$\gcapv{m}{i}{0.75}$:
%
\ylee{aes2.1a}
\gcapnI =
\gcapv{n-2d+2}{i_d}{2}\tcmp\cdots\tcmp
\gcapv{n-2}{i_2}{1.5}
\tcmp
\gcapv{n}{i_1}{0.75}.
\yeee
A \emph{\uptg} $\gcupnI$ is defined similarly:
%
\ylee{aes2.2a}
\gcupnI =
\gcupv{n}{i_1}{-0.75}
\tcmp
\gcupv{n-2}{i_2}{-1.25}
\tcmp
\cdots
\tcmp
\gcupv{n-2d+2}{i_d}{-2.25}
.
\yeee
The first proposition is obvious:
\begin{proposition}
\label{pr:3}
Every \TL\ \ttngnn\ $\xlam$ has a presentation
\xlee{aes2.3a}
\xlam = \gcupnIp\tcmp
\gcapnI,\qquad
\nI=\nIp.
\xeee
\end{proposition}
The number $\cpdlam=\nI=\nIp$ is determined by the tangle $\xlam$
and we call it  the \apdg\ (or \updg) of $\xlam$.

The second proposition is also obvious:
\begin{proposition}
\label{pr:1}
If at least one of two complexes $\xbCo$ and $\xbCt$ in $\dTLn$ is
\odct\ then their composition $\xbCo\tcmp\xbCt$ is \odct.
\end{proposition}
Note that even if
both complexes are \otbl,  then their composition is not necessarily
\otbl. Indeed, in contrast to the homological degree,
the \qdgr\ is not additive with respect to the composition of
tangles:
%
if  the composition of two \TL\ tangles contains a disjoint
circle then the \qdgr\ shifts of the rule\rx{ae1.01}
violate additivity. However, if the upper tangle in the composition has no caps or the
lower tangle has no cups then no circles are created and the
\otblc\ is maintained:
\begin{proposition}
\label{pr:2}
If a complex $\xbC$ in $\dTLv{n-2\cpdlam}$ is \otbl, then the complexes
$\ccupnI\tcmp\xbC$
and
$\xbC\tcmp\ccapnI$
are also \otbl.
\end{proposition}


\pr{Theorem}{th:2.1}
In order to construct the \odct\ and \otbl\ complex $\xbCmn$, we
use the presentation
\xlee{ae2.12}
\Cnv{\mrfm} \hteqv
\xbCon\tcmp\cbrmns\tgrsshomo,
\xeee
which follows from \eex{ae2.8} and\rx{ae2.6a}.
We construct $\xbCmn$ by
simplifying the complexes
$\bsymcats{\xlam\tcmp\gbrmn}$
for \TL\ \ttngnn s $\xlam$
appearing in the \qcmds\ of $\xbCon$, with the help of the
relation\rx{ae2.11}, thus creating necessary degree shifts, and then
using Corollary\rw{cr:otbl} which says that emerging \cbr s have \otbl\ categorification
complexes.




Let
$\dlam\tgrsshjmi$
be an object appearing in the $(-i)$-th \qcmd\ of
$\xbCon$ with a non-zero multiplicity (we made its homological degree explicit by
including $-i$ in the shift).
%
We apply \ex{ae2.11} consequently to every cap $\gcapnk$ appearing
in the \aptg\ $\gcapnI$ in the presentation\rx{aes2.3a} of
$\xlam$:
\begin{multline}
\label{ae2.13}
\dlam\tgrsshjmi\tcmp\cbrmns
\\
\hteqv
\Bigg(
\ccupnIp\tcmp
\bsymcatps{\gbrv{m}{n-2\cpdlam}{2.5}}\tcmp\ccapnI
\tgrsshv{\alm}{-\blm}^{2m} \tgrsshjmi\Bigg)
\tgrsshnontm,
\end{multline}
where
\xlee{ae2.14}
\alm = \sum_{k=1}^{\cpdlam-1}(n-2k)
,\qquad
\blm = \sum_{k=1}^{\cpdlam-1}(n-2k-1).
\xeee

The object $\dlam$ comes from the \odct\ complex $\xbCon$,
hence $\cpdlam>0$ and
the complex
in big brackets in the \rhs of \ex{ae2.13} is \odct\ in view of
Proposition\rw{pr:1}. Proposition\rw{pr:2} implies that the complex
$\ccupnIp\tcmp
\bsymcatps{\gbrv{m}{n-2\cpdlam}{2.5}}\tcmp\ccapnI$
is also \otbl. Since $\dlam$
comes from the \otbl\ complex $\xbCon$, the numbers $i$ and $j$
satisfy inequalities $i\geq 0$ and $i\leq j\leq 2i$. It is easy to check that the numbers
$\alm$ and $\blm$ of \ex{ae2.14} satisfy the same inequalities:
$\blm\geq 0$, $\blm\leq \alm \leq 2\blm$, hence the complex in big
brackets in the \rhs of \ex{ae2.13} is also \otbl. The
complex
$\xbCon\tcmp\cbrmns$
in the \rhs of \ex{ae2.12} is
composed of complexes\rx{ae2.13}, hence Theorem\rw{th:2.1} is proved.
\qed

\section{A categorified \JWp}
\label{s:prfs}

Consider the \chsq\rx{ae1.10c}. Theorem\rw{th:2.1} implies that
$\yordhr{\Cnv{\mrfm}} \geq 2m(n-1)+1$,
hence $\xctBn$ is
\Cch\ and it has a unique limit $\dlm\xctBn =\ctjwpn \in\dTLnp$.

Now we prove Theorems\rw{th:enum} and Theorem\rw{th:cnpr} which
describe the properties of $\ctjwpn$.

\pr{Theorem}{th:cnpr}
Consider the \chsq\rx{ae1.10c} truncated from below:
\ylee{eq:np1}
\xctBmn =
\Big(
\cbrmns \xraov{\mrfm}
\cbrmons \xrightarrow{\mrfmo}\cdots\Big)\longrightarrow\ctjwpn.
\yeee
According to Theorem\rw{th:rshfl}, the limit $\ctjwpn$ can be
presented as a cone\rx{ae2.m4}, where $\wbCmnp = \wbCmn\spshmnm$
and $\wbCmn$ is an infinite
\mtcn:
\begin{multline}
\nonumber
\wbCmn =\cdots\rightarrow\Cnv{\xbCvn{m+k}\tgrsshv{2kn}{-2k(n-1)+1}
\rightarrow
\cdots
\\
\cdots
\rightarrow
\Cnv{
\xbCvn{m+1}\tgrsshv{2n}{-2n+3}\rightarrow\xbCmn}
}
\end{multline}
with \odct\ and \otbl\ complexes $\xbCmn$ introduced in
Theorem\rw{th:2.1}. Hence the complex $\wbCmn$ itself is \odct\ and
\otbl.\qed

%

\pr{part 1 of Theorem}{th:enum}
The tangle composition with $\ccapni$ is a `continuous'
functor, that is, it can be applied to both sides of
$\dlm\xctBn = \ctjwpn$, hence $ \dlm\;
\ccapni\tcmp\xctBn = \ccapni\tcmp\ctjwpn$. According to \ex{ae2.11},
\begin{equation}
\nonumber
\begin{split}
\ccapni\tcmp\xctBn & = \Big(\ccapni\tcmp\cidbrn\rightarrow \cdots\rightarrow
\ccapni\tcmp\cbrmns\rightarrow\cdots\Big)
\\
& = \Big( \ccapni\rightarrow\cdots\rightarrow
\cbrmnmts\tcmp\ccapni
\spshmnm
\rightarrow\cdots
\Big).
\end{split}
\end{equation}
Since
\ylee{ae3.3}
\yordhb{\cbrmnmts\tcmp\ccapni
\spshmnm}
= 2m(n-1)\xrightarrow[m\rightarrow +\infty]{}+ \infty,
\yeee
according
to Theorem\rw{th:zl}, $\dlm\ccapni\tcmp\xctBn=0$, hence
$\ccapni\tcmp\ctjwpn$ is contractible.\qed

\begin{remark}
The contractibility of $\ctjwpn\tcmp\ccupni$ is proved similarly.
\end{remark}

\begin{corollary}
\label{cr:odct}
If $\xbC$ is a \odct\ complex in $\dTLnp$, then $\xbC\tcmp\ctjwpn$ is
contractible.
\end{corollary}

\pr{part 2 of Theorem}{th:enum}
According to
\ex{ae2.m5},
\begin{multline}
\nonumber
\ctjwpn\tcmp\ctjwpn \hteqv \CnBv{\wbCzn\qsho\longrightarrow\cidbrn}\tcmp\ctjwpn
\\
\hteqv\CnBv{\wbCzn\tcmp\ctjwpn\qsho\longrightarrow\cidbrn\tcmp\ctjwpn}
\hteqv\ctjwpn,
\end{multline}
where we used the fact that $\wbCzn$ is \odct\ and Corollary\rw{cr:odct} in order to establish the last
homotopy equivalence.\qed

\pr{Theorem}{th:alg}
The complexes $\ctjwpn$, $\wbCmn$ and $\cbrmns$ in \ex{ae2.m4} are
\otbl, hence they are \qpb\ and their $\Kz$ images are
well-defined. Applying $\Kz$ to this equation and taking into account \ex{eq:catKz} and
the definition\rx{ae1.10b1}, we find
\ylee{ae3.5}
\jwpn = q^{\vthf mn(n-1)}\abrmn - q^{2mn+1}\Kz(\wbCmn).
\yeee
The complex $\wbCmn$ is \otbl, so $\yordq{\Kz(\wbCmn)}\geq 0$ and
by Definition\rw{df:qlm} there is a limit\rw{ae1.9}.\qed

\section{The other projector}
A dual of an \ttngmn\ $\xtau$ is the \ttngnm\ tangle
$\xtaud$ which is its mirror image. Duality extends to an
isomorphism $\cTL \xrightarrow{\dsym} \cTLop$ combined with the
isomorphism of the ground ring $\Zqqi\xrightarrow{\dsym}\Zqqi$, such that
$\dsymv{q} = q^{-1}$. Furthermore, duality establishes an
isomorphism $\cTLpinf\xrightarrow{\dsym}\cTLminfop$, where
$\cTLminf$ is the analog of $\cTLpinf$ constructed over the ring
$\Zsqiq$ of Laurent series in $q^{-1}$.

Since the relations\rx{ae1.4} and\rx{ae1.4a} are dual to each other,
while the idempotency condition $\jwpn\tcmp\jwpn=\jwpn$ is duality
invariant, the uniqueness of the \JWp\ implies that it is duality
invariant: $\dsymv{\jwpn} = \jwpn$. Hence the corresponding
idempotents $\jwpnp\in\cTLpinf$ and $\jwpnm\in\cTLminf$ are also
dual to each other: $\jwpnm = \dsymv{(\jwpnp)}$. Taking the dual
of \ex{ae1.9} we find that $\jwpnm$ is the limit of \cbr s with
high positive (that is, \cclckw) twist:
\xlee{ae1.9b}
\lim_{m\rightarrow+\infty} q^{-\vthf mn(n-1)}\aobrmn = \jwpnm,
\xeee
because $\dsymv{\Big(\gbrmn\Big)} = \gobrmn$.


Duality extends to a contravariant
equivalence functor $\dTL\xrightarrow{\dsym}\dTLop$, where
$\dTLop$ is the same category as $\dTL$, except that the
composition of tangles is performed in reversed order. The functor
$\dsym$ also switches the signs of all three gradings of $\dTL$.
Applying the duality functor to the construction of $\ctjwpn$ we
find that there exists a \chsq
\begin{multline}
\label{ae1.10f}
\xctBnd =
\Big(
\cidbrn
\rxratv{\dmrfz}
\cobrons \rxratv{\dmrfo}
\cdots
\\
\cdots
\xleftarrow{\dmrfmmo}
\cobrmns \rxraov{\dmrfm}
\cobrmons \xleftarrow{\dmrfmo}\cdots\Big),
\end{multline}
where $-\xspsh$ denotes the grading shift which is opposite
to\rx{ae1.10b1}.
The system\rx{ae1.10f} is
dual to the system\rx{ae1.10c}
and it has an inductive limit $\ilm \xctBnd =\ctjwmn $, which satisfies
projector properties
\ylee{ae1.10f1}
\ctjwmn\tcmp \ctjwmn \hteqv \ctjwmn,\qquad
\ccapni \tcmp\ctjwmn \hteqv \ctjwmn\tcmp \ccupni\hteqv 0
\yeee
and has a presentation
\ylee{ae1.10f2}
\ctjwpn \hteqv \CnBv{ \dsymv{\wbCmn}\ospshmnm\qshmo
\xrahv{\dsymv{\chdlbtfm}} \cobrmns},
\yeee
%
where the complex $\wbCmn$ is
\odct\ and \otbl. In particular, at $m=0$ we get the dual of
presentation\rx{ae2.m5}:
\ylee{ae1.10f3}
\ctjwmn \hteqv \CnBv{ \dsymv{\wbCzn}\qshmo
\xrahv{\dsymv{\chdlbtfz}} \cidbrn},
\yeee
where the complex $\wbCzn$ is \odct\ and \otbl.





\end{document}

%% file: defs/gendef.tex
\def\ie{{\it i.e.\ }}

\def\cf{{\it cf.\ }}
\def\rhs{{\it r.h.s.\ }}

\def\deg{ \mathop{{\rm deg}}\nolimits }
\def\p{^{\prime}}

\def\pr#1#2{ \noindent{\em Proof of #1~\ref{#2}.} }
\def\proof{ \noindent{\em Proof.} }

\def\qed{ \hfill $\Box$ }

\def\lrbc#1{ \left( #1 \right) }
\def\lrbs#1{ \left[ #1 \right] }



%% file: defs/basdef.tex
\parskip=\medskipamount                 
\thicklines                     
\def\inbar{\vrule height1.5ex width.4pt depth0pt}
\def\IC{\relax\,\hbox{$\inbar\kern-.3em{\rm C}$}}
\def\IN{\relax{\rm I\kern-.18em N}}
\def\IQ{\relax\,\hbox{$\inbar\kern-.3em{\rm Q}$}}
\def\IR{\relax{\rm I\kern-.18em R}}
\def\ZZ{\relax{\sf Z\kern-.4em Z}}

\newtheorem{theorem}{Theorem}[section]
\newtheorem{proposition}[theorem]{Proposition}
\newtheorem{corollary}[theorem]{Corollary}

\newtheorem{lemma}[theorem]{Lemma}
\newtheorem{definition}[theorem]{Definition}
\newtheorem{remark}[theorem]{Remark}

\catcode`\@=11

\newif\if@fewtab\@fewtabtrue

\catcode`\@=11

\newif\if@fewtab\@fewtabtrue

{\count255=\time\divide\count255 by 60
\xdef\hourmin{\number\count255} \multiply\count255
by-60\advance\count255 by\time
\xdef\hourmin{\hourmin:\ifnum\count255<10 0\fi\the\count255}}
\def\ps@draft{\let\@mkboth\@gobbletwo
    \def\@oddhead{}
    \def\@oddfoot
      {\hbox to 7 cm{\footnotesize {\em Draft of \jobname:} \draftdate
       \hfil}\hskip -7cm\hfil\rm\thepage \hfil}
    \def\@evenhead{}\let\@evenfoot\@oddfoot}

\def\ceqno{\global\@fewtabfalse
    \ifcase\@eqcnt \def\@tempa{& & &}\or \def\@tempa{& &}
      \or \def\@tempa{&}
      \or\def\@tempa{}\fi\@tempa
{\rm(\theequation)}}

\def\aeqno#1{\global\@fewtabfalse
    \ifcase\@eqcnt \def\@tempa{& & &}\or \def\@tempa{& &}
      \or \def\@tempa{&}
      \or\def\@tempa{}\fi\@tempa
{\rm(\theequation,#1)}}

\def\label#1{\ifnum\draftcontrol=1
 \global\def\draftnote{$\scriptstyle #1$}\fi
 \@bsphack\if@filesw {\let\thepage\relax
   \def\protect{\noexpand\noexpand\noexpand}%
\xdef\@gtempa{\write\@auxout{\string
      \newlabel{#1}{{\@currentlabel}{\thepage}}}}}\@gtempa
   \if@nobreak \ifvmode\nobreak\fi\fi\fi
  \@esphack}

\def\alabel#1#2{\label{#1}\global\@fewtabfalse
    \ifcase\@eqcnt \def\@tempa{& & &}\or \def\@tempa{& &}
      \or \def\@tempa{&}
      \or\def\@tempa{}\fi\@tempa
{\hbox to 3cm{\phantom{\rm(\theequation,#2)} \draftnote
\hfil}\hskip -3cm {\rm(\theequation,#2)}}}

\def\clabel#1{\label{#1}\global\@fewtabfalse
    \ifcase\@eqcnt \def\@tempa{& & &}\or \def\@tempa{& &}
      \or \def\@tempa{&}
      \or\def\@tempa{}\fi\@tempa
{\hbox to 3cm{\phantom{\rm(\theequation)} \draftnote \hfil}\hskip
-3cm{\rm(\theequation)}}}

\def\eqnarray{\def\draftnote{{}}\global\@fewtabtrue
\stepcounter{equation}\let\@currentlabel=\theequation
\global\@eqnswtrue
\global\@eqcnt\z@\tabskip\@centering\let\\=\@eqncr
$$\halign to \displaywidth\bgroup\@eqnsel\hskip\@centering\@eqcnt\z@
  $\displaystyle\tabskip\z@{##}$&\global\@eqcnt\@ne
  \hskip 1\arraycolsep \hfil$\displaystyle{##}$\hfil
  &\global\@eqcnt\tw@ \hskip 1\arraycolsep
$\displaystyle\tabskip\z@{##}$ \hfil
\tabskip\@centering&\global\@eqcnt\thr@@\llap{##}\tabskip\z@ \cr}

\def\endeqnarray{\@@eqncr\egroup
      \global\advance\c@equation\m@ne$$\global\@ignoretrue}

\def\@eqnnum{\hbox to 3cm{\phantom{\rm(\theequation)} \draftnote
                         \hfil}\hskip -3cm {\rm(\theequation)}}

\def\@@eqncr{\let\@tempa\relax
    \ifcase\@eqcnt \def\@tempa{& & &}\or \def\@tempa{& &}
      \or \def\@tempa{&}
      \or\def\@tempa{}
\fi\@tempa \if@eqnsw \if@fewtab\@eqnnum\fi
\stepcounter{equation}\fi\global
\@eqnswtrue\global\@eqcnt\z@\global\@fewtabtrue\cr}

\def\draftcite#1{\ifnum\draftcontrol=1#1\else{}\fi}

\def\@lbibitem[#1]#2{\item{}\hskip -3cm \hbox to 2cm
{\hfil$\scriptstyle\draftcite{#2}$}\hskip
1cm[\@biblabel{#1}]\if@filesw
     {\def\protect##1{\string ##1\space}\immediate
      \write\@auxout{\string\bibcite{#2}{#1}}}\fi\ignorespaces}

\def\@bibitem#1{\item\hskip -3cm \hbox to 2cm
{\hfil $\scriptstyle\draftcite{#1}$}\hskip 1cm \if@filesw
\immediate\write\@auxout
       {\string\bibcite{#1}{\the\value{\@listctr}}}\fi\ignorespaces}

\def\draftdate{\number\month/\number\day/\number\year\ \ \ \hourmin }
 \global\def\draftcontrol{0}
\catcode`\@=12

\def\theequation{{\thesection.\arabic{equation}}}

\def\qq{\begin{eqnarray}}
\def\qqq{\end{eqnarray}}
\def\ee{\begin{eqnarray}}
\def\eee{\end{eqnarray}}
\def\rx#1{~(\ref{#1})}

\def\ex#1{eq.\hspace*{-3pt}\rx{#1}}
\def\eex#1{eqs.\hspace*{-3pt}\rx{#1}}
\def\cx#1{~\cite{#1}}
\def\rw#1{~\ref{#1}}

\def\xlee#1{ \begin{eqnarray} \label{#1} }
\def\xeee{ \end{eqnarray} }
\def\ylee#1{ \begin{eqnarray}\nonumber }
\def\yeee{ \end{eqnarray} }
\def\zlee#1{ \begin{displaymath} }
\def\zleee{ \end{displaymath} }
\def\wlee#1{ $ }
\def\weee{ $ }

\newlength{\shiftwidth}
\def\shift#1{&&\hbox to \shiftwidth{\hfill $\displaystyle#1$}}
\newlength{\sshiftwidth}
\def\sshift#1{\lefteqn{\hbox to
\sshiftwidth{\hfill$\displaystyle#1$}}}

\def\qbezier{\bezier{120}}
\def\DottedCircle{
\bezier{4}(0.966,-0.259)(1.04,0)(0.966,0.259)
\bezier{4}(0.966,0.259)(0.897,0.518)(0.707,0.707)
\bezier{4}(0.707,0.707)(0.518,0.897)(0.259,0.966)
\bezier{4}(0.259,0.966)(0,1.04)(-0.259,0.966)
\bezier{4}(-0.259,0.966)(-0.518,0.897)(-0.707,0.707)
\bezier{4}(-0.707,0.707)(-0.897,0.518)(-0.966,0.259)
\bezier{4}(-0.966,0.259)(-1.04,0)(-0.966,-0.259)
\bezier{4}(-0.966,-0.259)(-0.897,-0.518)(-0.707,-0.707)
\bezier{4}(-0.707,-0.707)(-0.518,-0.897)(-0.259,-0.966)
\bezier{4}(-0.259,-0.966)(0,-1.04)(0.259,-0.966)
\bezier{4}(0.259,-0.966)(0.518,-0.897)(0.707,-0.707)
\bezier{4}(0.707,-0.707)(0.897,-0.518)(0.966,-0.259) }
\def\Endpoint[#1]{
\ifcase#1 \put(1,0){\circle*{0.15}}
\or\put(0.866,0.5){\circle*{0.15}}
\or\put(0.5,0.866){\circle*{0.15}} \or\put(0,1){\circle*{0.15}}
\or\put(-0.5,0.866){\circle*{0.15}}
\or\put(-0.866,0.5){\circle*{0.15}} \or\put(-1,0){\circle*{0.15}}
\or\put(-0.866,-0.5){\circle*{0.15}}
\or\put(-0.5,-0.866){\circle*{0.15}} \or\put(0,-1){\circle*{0.15}}
\or\put(0.5,-0.866){\circle*{0.15}}
\or\put(0.866,-0.5){\circle*{0.15}} \fi}
\def\Arc[#1]{
\thicklines         
\ifcase#1 \bezier{25}(0.966,-0.259)(1.04,0)(0.966,0.259) \or
\bezier{25}(0.966,0.259)(0.897,0.518)(0.707,0.707) \or
\bezier{25}(0.707,0.707)(0.518,0.897)(0.259,0.966) \or
\bezier{25}(0.259,0.966)(0,1.04)(-0.259,0.966) \or
\bezier{25}(-0.259,0.966)(-0.518,0.897)(-0.707,0.707) \or
\bezier{25}(-0.707,0.707)(-0.897,0.518)(-0.966,0.259) \or
\bezier{25}(-0.966,0.259)(-1.04,0)(-0.966,-0.259) \or
\bezier{25}(-0.966,-0.259)(-0.897,-0.518)(-0.707,-0.707) \or
\bezier{25}(-0.707,-0.707)(-0.518,-0.897)(-0.259,-0.966) \or
\bezier{25}(-0.259,-0.966)(0,-1.04)(0.259,-0.966) \or
\bezier{25}(0.259,-0.966)(0.518,-0.897)(0.707,-0.707) \or
\bezier{25}(0.707,-0.707)(0.897,-0.518)(0.966,-0.259) \fi}
\def\DottedArc[#1]{
\ifcase#1 \bezier{4}(0.966,-0.259)(1.04,0)(0.966,0.259) \or
\bezier{4}(0.966,0.259)(0.897,0.518)(0.707,0.707) \or
\bezier{4}(0.707,0.707)(0.518,0.897)(0.259,0.966) \or
\bezier{4}(0.259,0.966)(0,1.04)(-0.259,0.966) \or
\bezier{4}(-0.259,0.966)(-0.518,0.897)(-0.707,0.707) \or
\bezier{4}(-0.707,0.707)(-0.897,0.518)(-0.966,0.259) \or
\bezier{4}(-0.966,0.259)(-1.04,0)(-0.966,-0.259) \or
\bezier{4}(-0.966,-0.259)(-0.897,-0.518)(-0.707,-0.707) \or
\bezier{4}(-0.707,-0.707)(-0.518,-0.897)(-0.259,-0.966) \or
\bezier{4}(-0.259,-0.966)(0,-1.04)(0.259,-0.966) \or
\bezier{4}(0.259,-0.966)(0.518,-0.897)(0.707,-0.707) \or
\bezier{4}(0.707,-0.707)(0.897,-0.518)(0.966,-0.259) \fi}
\def\Chord[#1,#2]{
\thinlines \ifnum#1>#2\Chord[#2,#1] \else\ifnum#1<#2 \ifcase#1
\ifcase#2 \or\qbezier(1,0)(0.516,0.138)(0.866,0.5)
\or\qbezier(1,0)(0.45,0.26)(0.5,0.866)
\or\qbezier(1,0)(0.327,0.327)(0,1)
\or\qbezier(1,0)(0.179,0.311)(-0.5,0.866)
\or\qbezier(1,0)(0.0536,0.2)(-0.866,0.5) \or\put(1, 0){\line(-2,
0){2}} \or\qbezier(1,0)(0.0536,-0.2)(-0.866,-0.5)
\or\qbezier(1,0)(0.179,-0.311)(-0.5,-0.866)
\or\qbezier(1,0)(0.327,-0.327)(0,-1)
\or\qbezier(1,0)(0.45,-0.26)(0.5,-0.866)
\or\qbezier(1,0)(0.516,-0.138)(0.866,-0.5) \fi \or\ifcase#2\or
\or\qbezier(0.866,0.5)(0.378,0.378)(0.5,0.866)
\or\qbezier(0.866,0.5)(0.26,0.45)(0,1)
\or\qbezier(0.866,0.5)(0.12,0.446)(-0.5,0.866)
\or\qbezier(0.866,0.5)(0,0.359)(-0.866,0.5)
\or\qbezier(0.866,0.5)(-0.0536,0.2)(-1,0) \or\put(0.866,
0.5){\line(-5, -3){1.73}}
\or\qbezier(0.866,0.5)(0.146,-0.146)(-0.5,-0.866)
\or\qbezier(0.866,0.5)(0.311,-0.179)(0,-1)
\or\qbezier(0.866,0.5)(0.446,-0.12)(0.5,-0.866)
\or\qbezier(0.866,0.5)(0.52,0)(0.866,-0.5) \fi \or\ifcase#2\or\or
\or\qbezier(0.5,0.866)(0.138,0.516)(0,1)
\or\qbezier(0.5,0.866)(0,0.52)(-0.5,0.866)
\or\qbezier(0.5,0.866)(-0.12,0.446)(-0.866,0.5)
\or\qbezier(0.5,0.866)(-0.179,0.311)(-1,0)
\or\qbezier(0.5,0.866)(-0.146,0.146)(-0.866,-0.5) \or\put(0.5,
0.866){\line(-3, -5){1}} \or\qbezier(0.5,0.866)(0.2,-0.0536)(0,-1)
\or\qbezier(0.5,0.866)(0.359,0)(0.5,-0.866)
\or\qbezier(0.5,0.866)(0.446,0.12)(0.866,-0.5) \fi
\or\ifcase#2\or\or\or \or\qbezier(0,1.)(-0.138,0.516)(-0.5,0.866)
\or\qbezier(0,1.)(-0.26,0.45)(-0.866,0.5)
\or\qbezier(0,1.)(-0.327,0.327)(-1,0)
\or\qbezier(0,1.)(-0.311,0.179)(-0.866,-0.5)
\or\qbezier(0,1.)(-0.2,0.0536)(-0.5,-0.866) \or\put(0, 1){\line(0,
-2){2}} \or\qbezier(0,1.)(0.2,0.0536)(0.5,-0.866)
\or\qbezier(0,1.)(0.311,0.179)(0.866,-0.5) \fi
\or\ifcase#2\or\or\or\or
\or\qbezier(-0.5,0.866)(-0.378,0.378)(-0.866,0.5)
\or\qbezier(-0.5,0.866)(-0.45,0.26)(-1,0)
\or\qbezier(-0.5,0.866)(-0.446,0.12)(-0.866,-0.5)
\or\qbezier(-0.5,0.866)(-0.359,0)(-0.5,-0.866)
\or\qbezier(-0.5,0.866)(-0.2,-0.0536)(0,-1) \or\put(-0.5,
0.866){\line(3, -5){1}}
\or\qbezier(-0.5,0.866)(0.146,0.146)(0.866,-0.5) \fi
\or\ifcase#2\or\or\or\or\or
\or\qbezier(-0.866,0.5)(-0.516,0.138)(-1,0)
\or\qbezier(-0.866,0.5)(-0.52,0)(-0.866,-0.5)
\or\qbezier(-0.866,0.5)(-0.446,-0.12)(-0.5,-0.866)
\or\qbezier(-0.866,0.5)(-0.311,-0.179)(0,-1)
\or\qbezier(-0.866,0.5)(-0.146,-0.146)(0.5,-0.866) \or\put(-0.866,
0.5){\line(5, -3){1.73}} \fi \or\ifcase#2\or\or\or\or\or\or
\or\qbezier(-1,0)(-0.516,-0.138)(-0.866,-0.5)
\or\qbezier(-1,0)(-0.45,-0.26)(-0.5,-0.866)
\or\qbezier(-1,0)(-0.327,-0.327)(0,-1)
\or\qbezier(-1,0)(-0.179,-0.311)(0.5,-0.866)
\or\qbezier(-1,0)(-0.0536,-0.2)(0.866,-0.5) \fi
\or\ifcase#2\or\or\or\or\or\or\or
\or\qbezier(-0.866,-0.5)(-0.378,-0.378)(-0.5,-0.866)
\or\qbezier(-0.866,-0.5)(-0.26,-0.45)(0,-1)
\or\qbezier(-0.866,-0.5)(-0.12,-0.446)(0.5,-0.866)
\or\qbezier(-0.866,-0.5)(0,-0.359)(0.866,-0.5) \fi
\or\ifcase#2\or\or\or\or\or\or\or\or
\or\qbezier(-0.5,-0.866)(-0.138,-0.516)(0,-1)
\or\qbezier(-0.5,-0.866)(0,-0.52)(0.5,-0.866)
\or\qbezier(-0.5,-0.866)(0.12,-0.446)(0.866,-0.5) \fi
\or\ifcase#2\or\or\or\or\or\or\or\or\or
\or\qbezier(0,-1.)(0.138,-0.516)(0.5,-0.866)
\or\qbezier(0,-1.)(0.26,-0.45)(0.866,-0.5) \fi
\or\ifcase#2\or\or\or\or\or\or\or\or\or\or
\or\qbezier(0.5,-0.866)(0.378,-0.378)(0.866,-0.5) \fi\fi\fi\fi}
\def\FullChord[#1,#2]{
\Endpoint[#1] \Endpoint[#2] \Arc[#1] \Arc[#2] \Chord[#1,#2] }
\def\EndChord[#1,#2]{
\Endpoint[#1] \Endpoint[#2] \Chord[#1,#2] }
\def\Picture#1{
\begin{picture}(2,1)(-1,-0.167)
#1
\end{picture}
}
\def\DottedChordDiagram[#1,#2]{
\Picture{\DottedCircle \FullChord[#1,#2]} }


%% file: defs/letdef.tex
\def\ZZ{ \mathbb{Z} }
\def\IQ{ \mathbb{Q} }
\def\IC{ \mathbb{C} }
\def\IR{ \mathbb{R} }


%% file: defs/locdef2m.tex
\def\tcmp{ \circ }

\def\plnr{planar}

\def\oltln{ \bigoplus_{\xlam\in\rTLn } }
\def\sltln{ \sum_{\xlam\in\rTLn } }

\def\lmii{ \lim_{i\rightarrow\infty} }

\def\xId{ \mathbbm{1} }

\def\JW{Jones-Wenzl}
\def\JWp{\JW\ projector}

\def\TLb{Temperley-Lieb}
\def\TLba{TL}
\def\TLa{\TLb\ algebra}

\def\cylb{torus}
\def\cbr{\cylb\ braid}

\def\aptg{cap-tangle}
\def\uptg{cup-tangle}

\def\dstt{distinguished triangle}

\def\stI{ \mathbf{I} }
\def\stIp{ \stI\p }

\def\ngbr{negative braid}

\def\nels#1{ |#1| }
\def\nI{ \nels{\stI} }
\def\nIp{ \nels{\stIp} }

\def\nsplcd{negatively spliced}

\def\apdg{cap-degree}
\def\updg{cup-degree}
\def\qdgr{$q$-degree}

\def\TNG{Tng}

\def\TL{TL}

\def\dTLt{ \widetilde{\dTL} }

\def\symalg#1{ \accentset{\smalg}{#1} }

\def\rxv#1{ \mathit{#1} }
\def\rTNG{ \rxv{\TNG} }

\def\rTL{ \rxv{\TL} }

\def\xop{ \mathrm{op} }

\def\cTL{ \symalg{\rTL} }
\def\cTLop{ \cTL^{\xop} }
\def\cTLv#1{ \cTL_{#1} }
\def\cTLn{ \cTLv{n} }
\def\cTLm{ \cTLv{m} }
\def\cTLvv#1#2{ \cTL\xivv{#1}{#2} }
\def\cTLmn{ \cTLvv{m}{n} }

\def\cTLnn{ \cTLvv{n}{n} }

\def\dTL{ \symcat{\rTL} }
\def\dTLop{ \dTL^{\xop} }
\def\dTLv#1{ \dTL_{#1} }
\def\dTLn{ \dTLv{n} }

\def\QcTL{ \mathrm{Q}\cTL }
\def\QcTLv#1{ \QcTL_{#1} }
\def\QcTLn{ \QcTLv{n} }

\def\cTLpinf{ \cTL^{+} }
\def\cTLminf{ \cTL^{-} }
\def\cTLminfop{ (\cTLminf)^{\xop} }

\def\symcat#1{ \accentset{\circ}{#1} }
\def\symbcat#1{ \lrbc{ #1 }^{\circ} }

\def\dTL{ \symcat{\rTL} }
\def\dTLnp{ \dTLp_n }

\def\dTLp{ \dTL^- }

\def\xivv#1#2{_{#1,#2} }
\def\xiv#1{_#1}

\def\rTNGvv#1#2{ \rTNG\xivv{#1}{#2} }
\def\rTNGv#1{ \rTNG\xiv{#1} }
\def\rTNGmn{ \rTNGvv{m}{n} }
\def\rTNGn{ \rTNGv{n} }

\def\rTLvv#1#2{ \rTL\xivv{#1}{#2} }
\def\rTLv#1{ \rTL\xiv{#1} }

\def\rTLn{ \rTLv{n} }

\def\ttngvv#1#2{$(#1,#2)$-tangle}
\def\ttngmn{\ttngvv{m}{n}}
\def\ttngnm{\ttngvv{n}{m}}
\def\ttngnn{\ttngvv{n}{n}}

\def\ttngzz{\ttngvv{0}{0}}

\def\Zqqi{ \ZZ[q,q^{-1}] }
\def\Zsqqi{ \ZZ[[q,q^{-1}] }
\def\Zsqiq{ \ZZ[[q^{-1},q] }
\def\Qq{ \mathrm{Q}(q) }

\def\hsq{ 1/2 }
\def\qh{ q^{\hsq} }
\def\qmh{ q^{-\hsq} }
\def\Zqqhi{ \ZZ[\qh,\qmh] }

\def\acapni{ \acapvv{n}{i} }

\def\ccapni{ \ccapvv{n}{i} }

\def\xU{ U }
\def\xUv#1#2{ \xU_{#1,#2} }
\def\xUni{ \xUv{n}{i} }

\def\jwp{ P }
\def\jwpv#1{ \jwp_{#1} }
\def\jwpn{ \jwpv{n} }
\def\jwpnp{ \jwpn^+ }
\def\jwpnm{ \jwpn^- }

\def\xtau{ T }

\def\xtauv#1{ \xtau_{#1} }
\def\xtauo{ \xtauv{1} }
\def\xtaut{ \xtauv{2} }

\def\xlam{ G }
\def\xlamp{ \xlam\p }
\def\ctau{ \symalg{\xtau} }
\def\clam{ \symalg{\xlam} }
\def\dtau{ \symcat{\xtau} }
\def\dlam{ \symcat{\xlam} }

\def\xnu{ \nu }

\def\xnuv#1{ \xnu_{#1} }
\def\xnum{ \xnuv{m} }

\def\xca{ a }
\def\xcav#1{ \xca_{#1} }
\def\xcal{ \xcav{\xlam} }
\def\xcalv#1{ \xcal(#1) }
\def\xcalt{ \xcalv{\xtau} }
\def\xcalw#1{ \xcav{\xlam,#1} }
\def\xcali{ \xcalw{i} }
\def\xcalj{ \xcalw{j} }
\def\xcaliv#1{ \xcali(#1) }
\def\xcalit{ \xcaliv{\xtau} }

\def\xL{ L }
\def\dL{ \symcat{\xL} }

\def\dLi{ \dLv{i} }

\def\dgo{ \deg_{\mathrm{h}} }
\def\dgt{ \deg_q }
\def\dgh{ \deg_{2} }

\def\tgrshv#1#2#3{ \lrbs{#2,#1,#3} }
\def\tgrsshv#1#2{ \lrbs{#2,#1} }

\def\tgrsshjmi{ \tgrsshv{j}{-i} }
\def\tgrsshomo{ \tgrsshv{1}{-1} }
\def\tgrsshnon{ \tgrsshv{n}{-n+1} }
\def\tgrsshnontm{ \tgrsshnon^{2m} }

\def\qshv#1{ \lrbs{ #1 }_{q} }

\def\qsho{ \qshv{1} }
\def\qshmo{ \qshv{-1} }

\def\spshnm{ \tgrsshv{n}{-n+1} }
\def\spshmnm{ \spshnm^{2m} }

\def\ospshnm{ \tgrsshv{n}{n-1} }
\def\ospshmnm{ \ospshnm^{2m} }

\def\Kz{\mathrm{K}_0}
\def\Kzg{$\Kz$-group}

\def\btcyl{ \brb_{\mathrm{cyl}} }

\def\btcylv#1{ \brb_{\mathrm{cyl},#1} }

\def\btcyln{ \btcylv{n} }

\def\vhf{ \frac{1}{2} }

\def\vthf{ \tfrac{1}{2} }
\def\vthh{ \tfrac{3}{2} }

\def\qpqi{ q + q^{-1} }
\def\qvh{ q^{\vhf} }
\def\qvmh{ q^{-\vhf} }

\def\xmrf{f}

\def\yal{ \alpha }
\def\ybet{ \beta }

\def\Opqv#1{ O_{+}(q^{#1})}

\def\Opqm{ \Opqv{m} }

\def\Oh{ O^{\mathrm{h}} }
\def\Ohp{ \Oh_{-} }

\def\Ohpv#1{ \Ohp(#1) }

\def\Ohpm{ \Ohpv{\ym} }
\def\ym{ m }
\def\ymv#1{ \ym_{#1} }
\def\ymi{ \ymv{i} }
\def\Ohpmi{ \Ohpv{\ymi} }

\def\symfr{ \blacksquare }
\def\symfr{ \circ }

\def\Cone{ \mathop{\mathrm{Cone}} }
\def\Cnv#1{ \Cone(#1) }
\def\CnBv#1{ \Cone\Big(#1\Big) }

\def\Cnbf{ \Cnv{\xbf} }
\def\Cnbfi{ \Cnv{\xbfi} }

\def\Cnbtfz{ \Cnv{\xbtfz} }
\def\Cnchdlbtfz{ \Cnv{\chdlbtfz} }

\def\xmu{ \mu }
\def\xcv#1{ c_{#1} }

\def\cjilam{ \xcv{i,j,\xmu}^{\xlam} }
\def\cjmilam{ \xcv{-i,j,\xmu}^{\xlam} }

\def\qpb{$q^+$\nobreakdash-bounded}

\def\ctjw{ \mathbf{P} }
\def\ctjwv#1{ \ctjw_{#1} }
\def\ctjwn{ \ctjwv{n} }
\def\ctjwp{ \ctjw^{-} }
\def\ctjwm{ \ctjw^{+} }
\def\ctjwpv#1{ \ctjwp_{#1} }
\def\ctjwmv#1{ \ctjwm_{#1} }
\def\ctjwpn{ \ctjwpv{n} }
\def\Kctjwpn{ \Kz(\ctjwpn) }
\def\ctjwmn{ \ctjwmv{n} }

\def\xsg{ S }
\def\xsgv#1{ \xsg_{#1} }
\def\xsgi{ \xsgv{i} }
\def\xsgiv#1{ \xsgv{i_{#1}} }
\def\xsgikpo{ \xsgiv{k+1} }

\def\brb{ B }
\def\brbp{ \brb\p }
\def\cbrb{ \symcat{\brb} }
\def\cbrba{ \symcats{\brb} }
\def\cbrbpa{ \symcats{\brbp} }
\def\spcc#1{ #1_{\sharp} }
\def\cbrbs{ \spcc{\cbrb} }
\def\cbrbas{ \spcc{\cbrba} }
\def\cbrbpas{ \spcc{\cbrbpa} }

\def\xC{ C }
\def\xCv#1{ \xC_{#1} }

\def\xCi{ \xCv{i} }
\def\xCmo{ \xCv{-1} }
\def\xCmt{ \xCv{-2} }
\def\xCmi{ \xCv{-i} }
\def\xCipo{ \xCv{i+1} }
\def\xbC{ \mathbf{C} }
\def\xbCv#1{ \xbC_{#1} }
\def\xbCo{ \xbCv{1} }
\def\xbCt{ \xbCv{2} }
\def\xbCvv#1#2{ \xbC_{#1,#2} }
\def\xbCvn#1{ \xbCvv{#1}{n} }
\def\xbCmn{ \xbCvn{m} }
\def\xbCon{ \xbCvn{1} }

\def\qmd{`module'}
\def\qmds{`modules'}
\def\qcmd{chain \qmd}
\def\qcmds{chain \qmds}

\def\otbl{angle-shaped}
\def\otblc{angle shape}

\def\dct{cut}
\def\dctv#1{#1-\dct}
\def\odct{\dctv{1}}

\def\crs{ n_{\times} }
\def\crsv#1{ \crs(#1) }

\def\mrf{ f }
\def\mrf{ \mathbf{f} }
\def\mrfv#1{ \mrf_{#1} }
\def\mrfz{ \mrfv{0} }
\def\mrfo{ \mrfv{1} }
\def\mrfm{ \mrfv{m} }
\def\mrfmmo{ \mrfv{m-1} }
\def\mrfmo{ \mrfv{m+1} }

\def\dmrf{ \dsymv{\mrf} }
\def\dmrfv#1{ \dmrf_{#1} }
\def\dmrfz{ \dmrfv{0} }
\def\dmrfo{ \dmrfv{1} }
\def\dmrfm{ \dmrfv{m} }
\def\dmrfmmo{ \dmrfv{m-1} }
\def\dmrfmo{ \dmrfv{m+1} }

\def\hteqv{ \sim }

\def\cpd{ d }
\def\cpdv#1{ \cpd_{#1} }
\def\cpdlam{ \cpdv{\xlam} }

\def\alm{ a_{\xlam} }
\def\blm{ b_{\xlam} }

\def\xnot{ \xlam_\varnothing }
\def\cnot{ \symcat{\xlam}_\varnothing }

\def\scA{ \mathcal{A} }
\def\scAs{ \scA_{\sharp} }

\def\scAp{ \scA\p }
\def\scAp{ \scAs }

\def\xctA{ \mathsf{A} }

\def\xbA{ \mathbf{A} }
\def\xbAs{ \xbA_{\sharp} }
\def\xA{ A }
\def\xAv#1{ \xA_{#1} }

\def\xAi{ \xAv{i} }
\def\xAio{ \xAv{i+1} }
\def\xAit{ \xAv{i+2} }

\def\xAmi{ \xAv{-i} }
\def\xAmio{ \xAv{-i+1} }

\def\xbAv#1{ \xbA_{#1} }
\def\xbAz{ \xbAv{0} }
\def\xbAo{ \xbAv{1} }
\def\xbAt{ \xbAv{2} }
\def\xbAi{ \xbAv{i} }
\def\xbAio{ \xbAv{i+1} }

\def\xbAp{ \xbA\p }
\def\xbApv#1{ \xbA_{\sharp,#1} }
\def\xbApz{ \xbApv{0} }
\def\xbApo{ \xbApv{1} }
\def\xbApt{ \xbApv{2} }
\def\xbApi{ \xbApv{i} }
\def\xbApio{ \xbApv{i+1} }

\def\xctB{ \mathcal{B} }
\def\xctBv#1{ \xctB_{#1} }
\def\xctBn{ \xctBv{n} }
\def\xctBnd{ \dsymv{\xctBn} }
\def\xctBmn{ \xctBv{m,n} }

\def\xbB{ \mathbf{B} }
\def\xB{ B }
\def\xBv#1{ \xB_{#1} }

\def\xBi{ \xBv{i} }
\def\xBio{ \xBv{i+1} }
\def\xBimo{ \xBv{i-1} }

\def\xd{ d }
\def\xdv#1{ \xd_{#1} }

\def\xdi{ \xdv{i} }
\def\xdmi{ \xdv{-i} }
\def\xdio{ \xdv{i+1} }

\def\xdimo{ \xdv{i-1} }

\def\xdp{ \xd\p }
\def\xdpv#1{ \xdp_{#1} }

\def\xdpi{ \xdpv{i} }
\def\xdpio{ \xdpv{i+1} }

\def\xdpimo{ \xdpv{i-1} }

\def\xbf{ \mathbf{f} }
\def\xbfv#1{ \xbf_{#1} }
\def\xbfi{ \xbfv{i} }

\def\xbfz{ \xbfv{0} }
\def\xbfo{ \xbfv{1} }
\def\xbft{ \xbfv{2} }

\def\xbh{ \mathbf{h} }
\def\xbhv#1{ \xbh_{#1} }
\def\xbhi{ \xbhv{i} }

\def\xbd{ \mathbf{d} }
\def\xbdv#1{ \xbd_{#1} }
\def\xbdi{ \xbdv{i} }

\def\xbhf{ \hat{\xbf} }
\def\xbhfv#1{ \xbhf_{#1} }
\def\xbhfi{ \xbhfv{i} }
\def\xbhfimo{ \xbhfv{i-1} }

\def\xbhh{ \hat{\xbh} }
\def\xbhhv#1{ \xbhh_{#1} }
\def\xbhhi{ \xbhhv{i} }
\def\xbhhz{ \xbhhv{0} }
\def\xbhho{ \xbhhv{1} }

\def\xbch{ \check{\xbh} }
\def\xbchv#1{ \xbch_{#1} }
\def\xbchi{ \xbchv{i} }

\def\xbhg{ \hat{\ybg} }
\def\xbhgv#1{ \xbhg_{#1} }
\def\xbhgi{ \xbhgv{i} }
\def\xbhgimo{ \xbhgv{i-1} }

\def\xbfAB{ \xbfv{\xbA\xbB} }
\def\xbfBC{ \xbfv{\xbB\xbC} }
\def\xbfAC{ \xbfv{\xbA\xbC} }

\def\xbtf{ \tilde{\xbf} }
\def\xbtfv#1{ \xbtf_{#1} }
\def\xbtfi{ \xbtfv{i} }
\def\xbtfio{ \xbtfv{i+1} }
\def\xbtfz{ \xbtfv{0} }

\def\xbtfm{ \xbtfv{m} }

\def\xbtfp{ \tilde{\xbf\p} }
\def\xbtfpv#1{ \xbtfp_{#1} }
\def\xbtfpi{ \xbtfpv{i} }

\def\xf{ f }
\def\xfv#1{ \xf_{#1} }
\def\xfi{ \xfv{i} }

\def\yf{ f }
\def\yfv#1{ \yf_{#1} }

\def\yfi{ \yfv{i} }
\def\yfio{ \yfv{i+1} }

\def\yg{ g }

\def\ybg{ \mathbf{\yg} }
\def\ybgv#1{ \ybg_{#1} }
\def\ybgA{ \ybgv{\xbA} }
\def\ybgB{ \ybgv{\xbB} }
\def\ybgi{ \ybgv{i} }
\def\ybgio{ \ybgv{i+1} }
\def\ybgimo{ \ybgv{i-1} }
\def\ybgz{ \ybgv{0} }
\def\ybgo{ \ybgv{1} }
\def\ybgt{ \ybgv{2} }
\def\ybgh{ \ybgv{3} }

\def\ybtg{ \tilde{\ybg} }
\def\ybtgv#1{ \ybtg_{#1} }
\def\ybtgi{ \ybtgv{i} }
\def\ybtgio{ \ybtgv{i+1} }

\def\ybtgz{ \ybtgv{0} }

\def\yh{ h }
\def\ybh{ \mathbf{\yh} }
\def\ybhv#1{ \ybh_{#1} }
\def\ybhi{ \ybhv{i} }
\def\ybhz{ \ybhv{0} }
\def\ybho{ \ybhv{1} }
\def\ybht{ \ybhv{2} }

\def\ybhp{ \ybh\p }
\def\ybhpv#1{ \ybhp_{#1} }
\def\ybhpi{ \ybhpv{i} }
\def\ybhpz{ \ybhpv{0} }
\def\ybhpo{ \ybhpv{1} }

\def\ybhtp{ \tilde{\ybh}\p }
\def\ybhtpv#1{ \ybhtp_{#1} }
\def\ybhtpi{ \ybhtpv{i} }

\def\xCh{ \mathop{\mathbf{Ch}} }
\def\xChv#1{ \xCh(#1) }
\def\xChA{ \xChv{\xctA} }

\def\xKh{ \mathop{\mathbf{K}} }
\def\xKhv#1{ \xKh(#1) }
\def\xKhA{ \xKhv{\xctA} }

\def\isor{isomorphism order}

\def\ysiov#1{ \left|\, #1\, \right|_{\cong} }
\def\ysiorv#1{ |\, #1\, |_{\cong} }
\def\ysiobf{ \ysiov{\xbf} }
\def\ysiobfi{ \ysiov{\xbfi} }

\def\stblz{stabilizing}

\def\xtrn{ \mathbf{t} }
\def\xtrnvv#1#2{ \xtrn^-_{\leq #1} #2 }
\def\xtrniv#1{ \xtrnvv{i}{#1} }
\def\xtrnNv#1{ \xtrnvv{N}{#1} }

\def\chcpl{chain complex}
\def\chmp{chain morphism}
\def\chcpls{\chcpl es}

\def\trnc{truncation}

\def\trncv#1{$#1$-th \trnc}
\def\trnci{\trncv{i}}

\def\chdl{ \delta }
\def\chdlv#1{ \chdl_{#1} }
\def\chdlbf{ \chdlv{\xbf} }
\def\chdlbfi{ \chdlv{\xbfi} }
\def\chdlbfz{ \chdlv{\xbfz} }
\def\chdlbtfz{ \chdlv{\xbtfz} }
\def\chdlbtfm{ \chdlv{\xbtfm} }

\def\idl{ \iota }
\def\idlv#1{ \idl_{#1} }
\def\idlbf{ \idlv{\xbf} }
\def\idlbgi{ \idlv{\ybgi} }
\def\idlbgz{ \idlv{\ybgz} }
\def\idlbgo{ \idlv{\ybgo} }

\def\idlbhi{ \idlv{\ybhi} }

\def\Cch{Cauchy}
\def\xsyst{system}
\def\xdsyst{direct \xsyst}
\def\Csq{\Cch\ \xsyst}

\def\dlm{ \lim\limits_{\longrightarrow} }
\def\ilm{ \lim\limits_{\longleftarrow} }

\def\qord{$q$-order}

\def\xsupv#1{ \sup\left\{#1 \right\} }
\def\xinfv#1{ \inf\left\{#1 \right\} }

\def\yordhb#1{ \bigl| \;#1 \;\bigr|_{\mathrm{h}} }
\def\yordh#1{ \left|\, #1\, \right|_{\mathrm{h}} }
\def\yordhr#1{ |\, #1\, |_{\mathrm{h} } }
\def\yordch#1{ \yordh{\Cnv{#1}} }
\def\yordq#1{ \left|\, #1\, \right|_q }

\def\zordh#1{ | #1 |_{\mathrm{h}} }
\def\zordch#1{ \zordh{\Cnv{#1}} }

\def\chsq{\xdsyst}

\def\chlm{ \lim\nolimits_{\xCh} }

\def\cmtr#1#2#3#4#5#6{
\xymatrix{
#4 \ar[r]_-{#1} \ar@/^1pc/[rr]^-{#3} &
#5 \ar[r]_-{#2} &
#6
},\qquad #3\hteqv #2\,#1
}

\def\atcmv#1#2{ #1\,#2 + #2\,#1 }

\def\xIdv#1{ \xId_{#1} }
\def\yIdAi{ \xIdv{\xbAi} }
\def\yIdA{ \xIdv{\xbA} }

\def\tchlm{chain limit}

\def\mmp{multi-map}

\def\ybC{ \mathbf{C} }
\def\ybCv#1{ \ybC_{#1} }
\def\ybCz{ \ybCv{0} }
\def\ybCo{ \ybCv{1} }
\def\ybCt{ \ybCv{2} }
\def\ybCi{ \ybCv{i} }
\def\ybCimo{ \ybCv{i-1} }

\def\ycC{ \mathcal{C} }
\def\yctC{ \tilde{\ycC} }

\def\ybtC{ \tilde{\ybC} }
\def\ybtCv#1{ \ybtC_{#1} }
\def\ybtCz{ \ybtCv{0} }
\def\ybtCo{ \ybtCv{1} }

\def\ybtCi{ \ybtCv{i} }

\def\ybtCio{ \ybtCv{i+1} }

\def\mtcn{multi-cone}

\def\wbC{ \tilde{\xbC} }
\def\wbCvv#1#2{ \wbC_{#1,#2} }
\def\wbCnv#1{ \wbCvv{#1}{n} }
\def\wbCmn{ \wbCnv{m} }
\def\wbCzn{ \wbCnv{0} }
\def\wbCmnp{ \wbCmn\p }

\def\xrahv#1{ \xrightarrow{\;\;\;#1\;\;\;} }
\def\xratv#1{ \xrightarrow{\;\;#1\;\;} }
\def\xraov#1{ \xrightarrow{\;#1\;} }

\def\rxratv#1{ \xleftarrow{\;\;#1\;\;} }
\def\rxraov#1{ \xleftarrow{\;#1\;} }

\def\Kh{Khovanov homology}

\def\clckw{clockwise}
\def\cclckw{counterclockwise}

\def\dsym{ \vee }
\def\dsymv#1{ #1^{\dsym} }
\def\xtaud{ \dsymv{\xtau} }

\def\pinft{\infty}

%% file: defs/locdef3.tex
\def\bldl{ \left\langle\!\left\langle }
\def\bldr{ \right\rangle\!\right\rangle }

\def\symcat#1{ \bldl #1 \bldr }
\def\symbcat#1{ \bldl #1 \bldr }
\def\Bsymcat#1{ \Bigl<\!\!\Bigl< #1 \Bigr>\!\!\Bigr> }
\def\bsymcat#1{ \bigl<\!\bigl< #1 \bigr>\!\bigr> }

\def\xspsh{ \mathrm{s} }
\def\ctss#1{ #1^{\xspsh } }
\def\ctsos#1{ #1^{-\xspsh } }
\def\symcats#1{ \ctss{\symcat{#1}} }

\def\bsymcats#1{ \ctss{\bsymcat{#1}} }

\def\bsymcatos#1{ \ctsos{\bsymcat{#1}} }

\def\ctpss#1{ \ctss{\spcc{#1}} }
\def\symcatps#1{ \ctpss{\symcat{#1}} }

\def\bsymcatps#1{ \ctpss{\bsymcat{#1}} }

\def\symalg#1{ \left\langle #1 \right\rangle }
\def\Bsymalg#1{ \Bigl< #1 \Bigr> }
\def\bsymalg#1{ \bigl< #1 \bigr> }

\def\dTL{ \mathsf{TL} }
\def\cTL{ \mathrm{TL} }
\def\rTL{ \mathit{TL} }

\def\acapni{ \acapvv{n}{i} }

\def\ccapni{ \ccapvv{n}{i} }

\def\acupni{ \acupvv{n}{i} }

\def\ccupni{ \ccupvv{n}{i} }

\def\xL{ L }
\def\dL{ \symcat{\xL} }

\def\dLi{ \dLv{i} }

\def\xnot{ \xlam_\varnothing }
\def\cnot{ \symcat{\xnot} }

\def\mpalg{ \symalg{-} }
\def\mpcat{ \symcat{-} }

\def\dtaus{ \symcats{\xtau} }

\def\dLi{ \xC_i }

\def\dLio{ \xC_{i+1} }

\def\xlam{ \lambda }
\def\xtau{ \tau }
\def\xsg{ \sigma }
\def\xbt{ \beta }
\def\brb{ \xbt }

%% file: defs/locdefgr3.tex
\def\lcir{
\xygraph{
!{0;/r1.5pc/:}
!{\vcap-}
!{\vcap}
}
}

\def\xcrsp{
\xygraph{
!{0;/r1.5pc/:}
[u(0.5)]
!{\xoverv}
}
}

\def\xpver{
\xygraph{
!{0;/r1.5pc/:}
[u(0.5)]
!{\xunoverv}
}
}

\def\xphor{
\xygraph{
!{0;/r1.5pc/:}
[u(0.5)]
!{\xunoverh}
}
}

\def\xvfro{
\xygraph{
!{0;/r1.5pc/:}
[u(0.5)]
!{\xcapv@(0)}
[u(1.5)]
[u(0.45)r(0.23)]
[d(1.5)]
*{\symfr\;\scriptstyle{1}}
}
}

\def\xvert{
\xygraph{
!{0;/r1.5pc/:}
[u(0.5)]
!{\xcapv@(0)}
}
}

\def\gbrv#1#2#3{
\xygraph{
!{0;/r1pc/:}
[u(0.5)]
!{\hcross}
!{\hcross}
[l(1)d(0.25)]
*{\scriptstyle{\vdots}}
[u(0.75)]
*{\scriptstyle{#1}}
[d(1)r(#3)]
*{\scriptstyle{#2}}
}
}

\def\gobrv#1#2#3{
\xygraph{
!{0;/r1pc/:}
[u(0.5)]
!{\hcrossneg}
!{\hcrossneg}
[l(1)d(0.25)]
*{\scriptstyle{\vdots}}
[u(0.75)]
*{\scriptstyle{#1}}
[d(1)r(#3)]
*{\scriptstyle{#2}}
}
}

\def\gbrmv#1#2{ \gbrv{m}{#1}{#2} }
\def\gbrmn{ \gbrmv{n}{1.5} }
\def\gbrmnmt{ \gbrmv{n-2}{2} }

\def\gobrmv#1#2{ \gobrv{m}{#1}{#2} }
\def\gobrmn{ \gobrmv{n}{1.5} }

\def\abrmn{ \bsymalg{\;\gbrmn} }
\def\cbrmns{ \bsymcats{\;\gbrmn} }
\def\cbrmnps{ \bsymcatps{\;\gbrmn} }

\def\aobrmn{ \bsymalg{\;\gobrmn} }
\def\cobrmns{ \bsymcatos{\;\gobrmn} }

\def\cbrmnmts{ \bsymcats{\;\gbrmnmt} }

\def\gbrov#1#2{ \gbrv{1}{#1}{#2} }
\def\gbron{ \gbrov{n}{1.5} }
\def\gbronmt{ \gbrov{n-2}{2} }

\def\gobrov#1#2{ \gobrv{1}{#1}{#2} }
\def\gobron{ \gobrov{n}{1.5} }

\def\cbrons{ \bsymcats{\;\gbron} }
\def\cbronps{ \bsymcatps{\;\gbron} }

\def\cobrons{ \bsymcatos{\;\gobron} }

\def\gbrmov#1#2{ \gbrv{m+1}{#1}{#2} }
\def\gbrmon{ \gbrmov{n}{1.5} }

\def\gobrmov#1#2{ \gobrv{m+1}{#1}{#2} }
\def\gobrmon{ \gobrmov{n}{1.5} }

\def\cbrmons{ \bsymcats{\;\gbrmon} }

\def\cobrmons{ \bsymcatos{\;\gobrmon} }

\def\gcapv#1#2#3{
\xygraph{
!{0;/r1pc/:}
[u(0.5)]
!{\hcap-}
[u(0.5)]
*{\scriptstyle{#2}}
[d(1)r(#3)]
*{\scriptstyle{#1}}
}
}

\def\gccapv#1#2#3{
\xygraph{
!{0;/r1pc/:}
[u(0.5)]
!{\hcap-}
[d(0.25)]
!{\hcap[-0.5]}
[u(0.25)]
[u(0.5)]
*{\scriptstyle{#2}}
[d(1)r(#3)]
*{\scriptstyle{#1}}
}
}

\def\gcapvi#1#2{ \gcapv{#1}{i}{#2} }
\def\gcapni{ \gcapvi{n}{0.5} }
\def\gcapnk{ \gcapv{n}{k}{0.5} }

\def\acapni{ \bsymalg{\gcapni} }
\def\ccapni{ \bsymcat{\!\gcapni} }

\def\gcapvI#1#2{ \gccapv{#1}{\stI}{#2} }
\def\gcapnI{ \gcapvI{n}{0.75} }
\def\ccapnI{ \bsymcat{ \gcapnI } }

\def\gcapfv#1#2#3#4{
\xygraph{
!{0;/r1pc/:}
[u(0.5)]
!{\hcap-}
[u(0.5)]
*{\scriptstyle{#2}}
[l(0.5)d(1)]
*{\circ}
[l(0.75)]
*{\scriptstyle{#3}}
[r(#4)]
*{\scriptstyle{#1}}
}
}

\def\gcapnit{ \gcapfv{n}{i}{2}{2} }
\def\gcapnik{ \gcapfv{n}{i}{k}{2} }

\def\gcupv#1#2#3{
\xygraph{
!{0;/r1pc/:}
[u(0.5)]
!{\hcap}
[u(0.5)]
*{\scriptstyle{#2}}
[d(1)r(#3)]
*{\scriptstyle{#1}}
}
}

\def\gccupv#1#2#3{
\xygraph{
!{0;/r1pc/:}
[u(0.5)]
!{\hcap}
[d(0.25)]
!{\hcap[0.5]}
[u(0.25)]
[u(0.5)]
*{\scriptstyle{#2}}
[d(1)r(#3)]
*{\scriptstyle{#1}}
}
}

\def\gcupvi#1#2{ \gcupv{#1}{i}{#2} }
\def\gcupni{ \gcupvi{n}{-0.5} }

\def\acupni{ \bsymalg{\gcupni} }
\def\ccupni{ \bsymcat{\gcupni\!} }

\def\gcupvI#1#2{ \gccupv{#1}{\stI}{#2} }
\def\gcupnI{ \gcupvI{n}{-0.5} }
\def\ccupnI{ \bsymcat{ \gcupnI } }

\def\gcupvIp#1#2{ \gccupv{#1}{\stI\p}{#2} }
\def\gcupnIp{ \gcupvIp{n}{-0.5} }
\def\ccupnIp{ \bsymcat{\gcupnIp} }

\def\gidbrv#1#2{
\xygraph{
!{0;/r1pc/:}
[u(0.8)]
!{\xcaph@(0)}
[d(1.5)l(1)]
!{\xcaph@(0)}
[u(1)l(0.4)]
*{\vdots}
[d(0.3)r(#2)]
*{\scriptstyle{#1}}
}
}

\def\gidbrn{ \gidbrv{n}{1} }

\def\cidbrn{ \bsymcat{\,\gidbrn} }

\def\tstgrt{
\xygraph{
!{0;/r1pc/:}
!{\xcaph@(0)}
[d(1.5)l(1)]
!{\xcaph@(0)}
[u(1)l(0.4)]
*{\vdots}
[d(0.3)r(1)]
*{\scriptstyle{n}}
}
}

\def\testgr{
\underbrace{
\xygraph{
!{0;/r1.5pc/:}
[] *+[F]{\;\;\scriptstyle{n}\;\;}
[l(0.5)d(0.4)]
!{\xcapv[0.5]@(0)}
[r(1)u(1)]
!{\xcapv[0.5]@(0)}
[u(2.3)]
!{\xcapv[0.5]@(0)}
[l(1)u(1)]
!{\xcapv[0.5]@(0)}
}
}_{n}
}

\def\testgr1{
\xygraph{
!{0;/r1.5pc/:}
!{\vtwist}
!{\vtwist}
}
}

\def\testgr2{
\xygraph{
!{0;/r1pc/:}
!{\hcrossneg}
!{\hcrossneg}
}
}